\documentstyle[12pt,amssymb,amscd]{amsart}
\newtheorem{thm}{Theorem}[section]
\newtheorem{dfn}{Definition}[section]
\newtheorem{lem}{Lemma}[section]
\newtheorem{cor}{Corollary}[section]
\newtheorem{prop}{Proposition}[section]

\newcommand{\xt}{\mbox{$\tilde{X}$}}
\newcommand{\yt}{\mbox{$\tilde{Y}$}}
\newcommand{\dt}{\mbox{$\tilde{D}$}}
\newcommand{\ct}{\mbox{$\tilde{C}$}}

\begin{document}


\title{General $n$-canonical divisors on two-dimensional smoothable
semi-log-terminal singularities}
\author{Masayuki Iwamoto}
\maketitle
\section{Introduction}
 This paper is devoted to some fundumental calculation on
2-dimensional smoothable semi-log-terminal singularities.
If we study minimal or canonical models of one parameter degeneration
of algebraic surfaces, we must treat singularities that appear in
the central fiber.
Smoothable semi-log-terminal singularities are the singularities
of the central fiber of the minimal model of degeneration,
and the singularities of the central fiber of the canonical model of
degeneration
which may have large Gorenstein index.
Koll\'{a}r and Shepherd-Barron caracterized these singularities
in \cite{ksb}, but for numerical theory of degeneration,
we need more detailed information.
\par
In this paper, we calculate general $n$-canonical divisors on
these singularities, in other words, we calculate the full sheaves
associated to the double dual of the $n$-th tensor power of the
dualizing sheaves. And the application of this result,
we bound the Gorenstein index by the local self intersection number
of the $n$-canonical divisor.
\par
Notation: In this paper,
\begin{align*}
[q_1, q_2, q_3, \dots]&=q_1+
\cfrac{1}{q_2+
\cfrac{1}{q_3+{}_{\ddots}
}}\\
[[q_1, q_2, q_3, \dots]]&=q_1-
\cfrac{1}{q_2-
\cfrac{1}{q_3-{}_{\ddots}
}}
\end{align*}
If $p:\tilde{X}\to X$ is a birational morphism and $D$ is a
divisor on $X$,
we denote by $\tilde{D}$ the proper transform of $D$ on $\tilde{X}$.

\section{Basic calculation}
Let $Y$ be a cyclic quotient singularity of the form Spec$(\Bbb
C[z_1,z_2])/\langle\alpha\rangle$,
where $\langle\alpha\rangle$ is a cyclic group of order $r$ and $\alpha$ acts
on Spec$(\Bbb C[z_1,z_2])/\langle\alpha\rangle$ as $(\alpha^* z_1,\alpha^*
z_2)=(\eta^s z_1,\eta z_2)$
in which $\eta$ is a primitive $r$-th root of unity, $(r,s)=1$,
and $0<s<r$.
Let $\dfrac{r}{s}=[[q_1, q_2, \dots, q_k]]$ be an expansion into continued
fraction, and $r_{i}$ be the $i$-th remainder of the Euclidean algorithm, i.e.
$\{r_i\}_{i=0, 1, \dots, k+1}$ is a seqence determined by
$r_0=r, r_1=s, r_{i-1}=q_i r_i-r_{i+1}$.
Let $\dfrac{P_i}{Q_i}$ be the $i$-the convergent,
i.e. $\{P_i\}_{i=-1, 0, \dots, k}$
 is a sequence determined by $P_{-1}=0, P_0=1, P_i=q_i P_{i-1}-P_{i-2}$ and
$\{Q_i\}_{i=-1, 0, \dots, k}$
 is determined by $Q_{-1}=-1, Q_0=0, Q_i=q_iQ_{i-1}-Q_{i-2}$.
Let $f:\Bbb C^{2}\rightarrow Y$ be a quotient map, and $p:\yt\rightarrow Y$
 the minimal desingularization.
It is well known that the dual graph of the exceptional divisors of $p$
 is a chain of rational curves $\cup_{1\le i\le k}E_i$ such that
$E_i^2=-q_i$.
We put
$$\lambda_i=P_{i-1},\quad  \mu_{i}=r_i,\quad  \alpha^j_i=
\begin{cases}
\dfrac{1}{m}\mu_i\lambda_j\ &(\text{for}\ j\le i)\\
\dfrac{1}{m}\lambda_i\mu_j\ &(\text{for}\ i<j)
\end{cases}
$$
Note that $c_1z_1^{\lambda_i}+c_2z_2^{\mu_i}$
is a $\langle\alpha\rangle$-semi-invariant since $r_1P_i-r_0Q_i=r_{i+1}$.

\begin{lem}\label{lemcurve}
Let $C_i$ be a divisor on $Y$ such that
$f^*C_i=(c_iz_1^{\lambda_i}+c_2z_2^{\mu_i}=0)$
 in which $c_1,c_2\in \Bbb C^*$.
Then
$$
\text{\em{(i)}}\ \ct_i\cdot E_j=\delta_{i,j}\quad
\text{\em{(ii)}}\ p^*C_i=\ct_i+\sum\alpha^i_jE_j$$
\end{lem}
\begin{pf}
 If we write $Y=T_N\text{emb}(\sigma)$, where $N=\Bbb Z n_1+\Bbb Z n_2$
and $\sigma=\Bbb R_{\geq 0}n_1+\Bbb R_{\geq 0}[(r-s)n_1+n_2]$,
then $E_i$ corresponds to $(P_{i-1}-Q_{i-1})n_1+P_{i-1}n_2$.
The rest of proof is a direct calculation using the above description and
the formula $P_iQ_{i-1}-Q_iP_{i-1}=-1$, and it can be easily done.
\end{pf}

 Next lemma is a easy fact on continued fraction.
We denote by $(2, q)$ the sequence $(2, 2, \dots, 2)$ of length $q$.
\begin{lem}\label{lemcont1}
Let $a$, $m$ be a natural number such that $(a, m)=1$,
$\frac{m}{2}<a<m$. Put $b=m-a$.
Let $\dfrac{a}{b}=[q_1, q_2, \dots, q_k]$.
Then
\begin{enumerate}
\item If $k$ is even,
\begin{align*}
\dfrac{m}{a}=&[[(2, q_1), q_2+2, (2, q_3-1), q_4+2, \dots, (2, q_{k-3}-1),
q_{k-2}+2,\\
& (2, q_{k-1}-1), q_k+1]]\\
\dfrac{m}{b}=&[[q_1+2, (2, q_2-1), q_3+2, (2, q_4-1), q_5+2, \dots,
(2, q_{k-2}-1),\\
& q_{k-1}+2, (2, q_k+1)]]
\end{align*}
\item If $k$ is odd,
\begin{align*}
\dfrac{m}{a}=&[[(2, q_1), q_2+2, (2, q_3-1), q_4+2, \dots, q_{k-3}+2,
(2, q_{k-2}-1),\\
&q_{k-1}+2, (2, q_k-1)]]\\
\dfrac{m}{b}=&[[q_1+2, (2, q_2-1), q_3+2, (2, q_4-1), q_5+2, \dots,
(2, q_{k-3}-1), q_{k-2}+2,\\
& (2, q_{k-1}-1), q_k+1]]
\end{align*}
\end{enumerate}
\end{lem}
\begin{pf}
Since this is elementary we left it for the reader.
\end{pf}

 In the rest of this section we shall index the exceptional divisors
of the minimal (semi-) resolution smoothable of
the semi-log-terminal singularity.
\par
First we treat the normal case.
Let $(a, d, m)$ be a triplet of positive integers such that $a<b$ and
$a$ is prime to $m$.
We denote by $X_{a, d, m}$ a 2-dimentional quotient singularity
of the form Spec$(\Bbb C[z_1, z_2])/\langle\alpha\rangle$,
where $\langle\alpha\rangle$ is a cyclic group of order $dm^2$ and
$\alpha$ acts on Spec$(\Bbb C[z_1, z_2])$ as
$\alpha^*(z_1, z_2)=(\varepsilon^{adm-1}z_1, \varepsilon z_2)$
in which $\varepsilon$ is a primitive $dm^2$-th root of unity.
By [K-SB Proposition 3.10], a singularity of class T which is not
RDP is analytically isomorphic to $X_{a, d, m}$ for
some $(a, d, m)$.
Let $f:\Bbb C^2\to X_{a, d, m}$ be a quotient map and
$p:\xt_{a, d, m}\to X_{a, d, m}$
the minimal desingularization.
We assume $2a>m$ since $X_{a, d, m}\simeq X_{m-a, d, m}$.
Put $b=m-a$, and let $\dfrac{a}{b}=[q_1, q_2, \dots, q_k]$
be an expantion into continued fraction.
Let $r_i$ be the $i$-th remainder of the Euclidean algorithm, i.e.
 $\{r_i\}_{i=0, 1, \dots, k+1}$ is a sequence determined by
$r_0=a, r_1=m-a, r_{i-1}=q_i r_i+r_{i+1}$.
Let $\dfrac{P_i}{Q_i}$ be an the $i$-th convergent,
i.e. $\{P_i\}_{i=-1, 0, \dots, k}$ is a sequence
determined by $P_{-1}=0, P_0=1, P_i=q_iP_{i-1}+P_{i-2}$
and $\{Q_i\}_{i=-1, 0, \dots, k}$ is determined by
$Q_{-1}=1, Q_0=0, Q_i=q_iQ_{i-1}+Q_{i-2}$.
\begin{lem}\label{lemcont2}
Let $m$, $a$, $b$, $d$ be positive integers such that $m=a+b$,
 $a>b$, $(a, b)=1$.
Let $\dfrac{a}{b}=[q_1, q_2, \dots, q_k]$ be an expansion into
continued fraction.
Then $$\dfrac{dma-1}{dmb+1}=[q'_1, q'_2, \dots, q'_{k'}]$$
where $q'_i$ is as follows.
\begin{enumerate}
\item If $d=1$ and $k$ is even,
\begin{equation*}
q'_i=
\begin{cases}
q_1\quad &(i=1)\\
q_i\quad &(2\le i\le k-1)\\
q_k+1\quad &(i=k)\\
q_k-1\quad &(i=k+1)\\
q_{2k-i+1}\quad &(k+2\le i\le 2k-1)\\
q_1+1\quad &(i=2k=k')
\end{cases}
\end{equation*}
\item If $d=1$ and $k$ is odd,
\begin{equation*}
q'_i=
\begin{cases}
q_1\quad &(i=1)\\
q_i\quad &(2\le i\le k-1)\\
q_k-1\quad &(i=k)\\
q_k+1\quad &(i=k+1)\\
q_{2k-i+1}\quad &(k+2\le i\le 2k-1)\\
q_1+1\quad &(i=2k=k')
\end{cases}
\end{equation*}
\item If $d\ge 2$ and $k$ is even,
\begin{equation*}
q'_i=
\begin{cases}
q_1\quad &(i=1)\\
q_i\quad &(2\le i\le k)\\
d-1\quad &(i=k+1)\\
1\quad &(i=k+2)\\
q_k-1\quad &(i=k+3)\\
q_{2k+3-i}\quad &(k+4\le i\le 2k+1)\\
q_1+1\quad &(i=2k+2=k')
\end{cases}
\end{equation*}
\item If $d\ge 2$ and $k$ is odd,
\begin{equation*}
q'_i=
\begin{cases}
q_1\quad &(i=1)\\
q_i\quad &(2\le i\le k-1)\\
q_k-1\quad &(i=k)\\
1\quad &(i=k+1)\\
d-1\quad &(i=k+2)\\
q_{2k+3-i}\quad &(k+3\le i\le 2k+1)\\
q_1+1\quad &(i=2k+2=k')
\end{cases}
\end{equation*}
\end{enumerate}
\end{lem}
By Lemma \ref{lemcont1} and \ref{lemcont2}, we can calculate the dual graph of
 the exceptional divisors of $p$ in terms of the continued
fraction expansion of $\dfrac{a}{b}$.
(See [K-SB])
{}From now, we assume $k$ is even since the calculation is the same for
odd $k$.
We shall index the exceptional divisors in the following manner.
Set the index set $I_o, I_e, I$ as follows:
\begin{align*}
I_o&=\{(i, j)|1\leq i\leq k+1;i\ \text{odd};
1\leq j\leq q_i(\text{for}\ i<k+1),
1\leq j\leq d(\text{for}\ i=k+1)\}\\
I_e&=\{(i, j)|1\leq i\leq k+1;i\ \text{even};
1\leq j\leq q_i(\text{for}\ i<k),
1\leq j\leq q_k-1(\text{for}\ i=k)\}\\
I&=I_o\amalg I_e
\end{align*}
We define $\rho^i_j$, $\bar{\lambda}^i_j$, $\hat{\lambda}^i_j$,
$\lambda^i_j$, $\bar{\mu}^i_j$, $\hat{\mu}^i_j$, $\mu^i_j$
for $(i, j)\in I$ as follows.
\begin{equation*}
\rho^i_j=r_{i-1}-(j-1)r_i
\end{equation*}
\begin{equation*}
\bar{\lambda}^i_j=
\begin{cases}
P_{i-2}+(j-1)P_{i-1}\quad &((i, j)\in I_o)\\
-\{P_{i-2}+(j-1)P_{i-1}\}+da\rho^i_j\quad &((i, j)\in I_e)
\end{cases}
\end{equation*}
\begin{equation*}
\hat{\lambda}^i_j=
\begin{cases}
Q_{i-2}+(j-1)Q_{i-1}\quad &((i, j)\in I_o)\\
-\{Q_{i-2}+(j-1)Q_{i-1}\}+db\rho^i_j\quad &((i, j)\in I_e)
\end{cases}
\end{equation*}
\begin{equation*}
\bar{\mu}^i_j=
\begin{cases}
-\{P_{i-2}+(j-1)P_{i-1}\}+da\rho^i_j\quad &((i, j)\in I_o)\\
P_{i-2}+(j-1)P_{i-1}\quad &((i, j)\in I_e)
\end{cases}
\end{equation*}
\begin{equation*}
\hat{\mu}^i_j=
\begin{cases}
-\{Q_{i-2}+(j-1)Q_{i-1}\}+db\rho^i_j\quad &((i, j)\in I_o)\\
Q_{i-2}+(j-1)Q_{i-1}\quad &((i, j)\in I_e)
\end{cases}
\end{equation*}
\begin{equation*}
\lambda^i_j=\bar{\lambda}^i_j+\hat{\lambda}^i_j,\quad
\mu^i_j=\bar{\mu}^i_j+\hat{\mu}^i_j
\end{equation*}
We can write $Y=T_N\text{emb}(\sigma)$, where
$$N=\Bbb Z n_1+\Bbb Z n_2,\quad \sigma=\Bbb R_{\ge 0}n_1+\Bbb R_{\ge
0}[d(bm+1)n_1+dm^2n_2]$$
We denote by $E^i_j$ the exceptional divisor associated to
$\hat{\lambda}^i_jn_1+\lambda^i_jn_2$.
Note that by Lemma \ref{lemcurve}, for $\iota\in I$,
the proper transform of
$C_{\iota}=(c_1z_1^{\lambda_{\iota}}+c_2z_2^{\mu_{\iota}})/
\langle\alpha\rangle\in X_{a,d,m}$
intersects the exceptional locus transversely at $E_{\iota}$.
We define the order `$\le$' in the index set $I$ by
the lexicographic order.
\par

 Next we treat the non-normal case.
 We denote by $NC^2=\text{Spec}\Bbb C [z_1, z_2, z_3]/(z_1 z_2)$
 a 2-dimentional normal crossing point.
Let $(a, m)$ be a pair of positive integers
such that $0<a<m$ and $a$ is prime to $m$.
Put $b=m-a$, and let $a'$ and $b'$ be integers such that
$aa'\equiv bb'\equiv 1\pmod{m}$, $0<a'<m$, and $0<b'<m$.
Let $\langle \alpha \rangle$ be a cyclic group of order $m$,
and let $\langle \alpha \rangle$ act on $NC^2$ as
$(\alpha^*z_1, \alpha^*z_2, \alpha^*z_3)
=(\varepsilon^{a'}z_1, \varepsilon^{b'}z_2, \varepsilon z_3)$
where $\varepsilon$ is a primitive $m$-th root of unity.
We denote by $X_{a, m}$ the quotient of $NC^2$ by
this $\langle \alpha \rangle$-action.
By [K-SB], 2-dimentional smoothable semi-log-terminal singularity
which is neither normal nor $NC^2$ is analytically isomorphic to
$X_{a, m}$ for some $(a, m)$. Put $X=X_{a, m}$.
Let $f:NC^2\to X$ be the quotient map and $p:\tilde{X}\to X$
the minimal semi-resolution.
Let $g:\Bbb C^2_o\amalg\Bbb C^2_e\to NC^2$,
$g_X:X_o\amalg X_e\to X$,
and $g_{\tilde{X}}:\tilde{X}_o\amalg\tilde{X}_e\to\tilde{X}$
be normalizations,
where $\Bbb C^2_o=\text{Spec}\Bbb C [z_1, z_3]$,
$\Bbb C^2_e=\text{Spec}\Bbb C [z_2, z_3]$,
$X_o$ (resp. $X_e$) is the quotient of $\Bbb C^2_o$ (resp. $\Bbb C^2_e$),
 and $\tilde{X}_o$ (resp. $\tilde{X}_e$) is the minimal resolution
of $X_o$ (resp. $X_e$). We get the following diagram.

\begin{equation*}
 \begin{CD}
 \Bbb C^2_o\amalg \Bbb C^2_e @>{f_o\amalg f_e}>>
 X_o\amalg X_e @<{p_o\amalg p_e}<<
 \tilde{X}_o\amalg\tilde{X}_e \\
 @V{g}VV @VV{g_X}V @VV{g_{\tilde{X}}}V \\
 NC^2 @>>{f}> X @<<{p}< \tilde{X}
 \end{CD}
\end{equation*}

We denote by $\Delta$, $\Delta '$, and
$\tilde{\Delta}$ the double curve of
$X$, $NC^2$, and $\tilde{X}$ respectively.
Let $\Delta_o$ (resp. $\Delta_e$) be the inverse image of
$\Delta$ in $X_o$ (resp. $X_e$), and define
$\Delta '_o$, $\Delta '_e$, $\tilde{\Delta}_o$, and
$\tilde{\Delta}_e$ similarly.
 We assume $k$ is even.
 Set the index set $I_o$, $I_e$, $I$ as follows:
\begin{align*}
I_o&=\{(i, j)|1\le i\le k+1;i\ \text{odd};
1\le j\le q_i(\text{for}\ i<k), j=1(\text{for}\ i=k+1)\}\\
I_e&=\{(i, j)|1\le i\le k+1;i\ \text{even};
1\le j\le q_i\}\\
I&=I_o\amalg I_e
\end{align*}
We define $\lambda^i_j$, $\mu^i_j$ for $(i, j)\in I$ as follows:
\begin{equation*}
\lambda^i_j=
\begin{cases}
P_{i-2}+Q_{i-2}+(j-1)(P_{i-1}+Q_{i-1})\quad &((i, j)\in I_o)\\
r_{i-1}-(j-1)r_i\quad &((i, j)\in I_e)
\end{cases}
\end{equation*}
\begin{equation*}
\mu^i_j=
\begin{cases}
r_{i-1}-(j-1)r_i\quad &((i, j)\in I_o)\\
P_{i-2}+Q_{i-2}+(j-1)(P_{i-1}+Q_{i-1})\quad &((i, j)\in I_e)
\end{cases}
\end{equation*}
We define $\bar{\lambda}^i_j$, $\hat{\lambda}^i_j$ for $(i, j)\in I_o$
as follows:
$$\bar{\lambda}^i_j=P_{i-2}+(j-1)P_{i-1},\quad
\hat{\lambda}^i_j=Q_{i-2}+(j-1)Q_{i-1}$$
We define $\bar{\mu}^i_j$, $\hat{\mu}^i_j$ for $(i, j)\in I_e$ as follows:
$$\bar{\mu}^i_j=P_{i-2}+(j-1)P_{i-1},\quad
\hat{\mu}^i_j=Q_{i-2}+(j-1)Q_{i-1}$$
We write $X_o=T_{N^o}\text{emb}(\sigma^o)$ and
 $X_e=T_{N^e}\text{emb}(\sigma^e)$ where
$$N^o=\Bbb Z n^o_1+\Bbb Z n^o_2,\quad
\sigma^o=\Bbb R_{\ge 0}n^o_1+\Bbb R_{\ge 0}(bn^o_1+mn^o_2)$$
$$N^e=\Bbb Z n^e_1+\Bbb Z n^e_2,\quad
\sigma^e=\Bbb R_{\ge 0}n^e_1+\Bbb R_{\ge 0}(an^e_1+mn^e_2)$$
For $(i, j)\in I_o$ (resp. $\in I_e$), we denote by $E^i_j$
the exceptional divisor of $p_o$ (resp. $p_e$) which
associated to $\hat{\lambda}^i_j n^o_1+\lambda^i_j n^o_2\in N^o$
(resp. $\hat{\mu}^i_j n^e_1+\mu^i_j n^e_2\in N^e$).
Note that by Lemma \ref{lemcurve},
for $\iota\in I_o$ (resp.$I_e$),
the proper transform of
$C_{\iota}=(c_1z_3^{\lambda_{\iota}}+c_2z_1^{\mu_{\iota}})/
\langle\alpha\rangle\in X_o$ (resp.
$C_{\iota}=(c_1z_2^{\lambda_{\iota}}+c_2z_3^{\mu_{\iota}})/
\langle\alpha\rangle\in X_e$)
intersects the exceptional locus transversly at $E_{\iota}$.
We define the order in $I$ as the same way as the normal case.

\par
In the rest of this paper, we treat $X_{a,d,m}$ and $X_{a,m}$
simultaneously, otherwise we specifically state the normal or
non-normal case.

\section{$\lambda$-expansion and $\mu$-expantion}

In this section we introduce the notion of
$\lambda$-expansion and $\mu$-expansion, which is the key
in this paper.

\begin{dfn}\label{dfnlambda}
Let $L=(j_1, l_2, j_3, l_4, \dots, j_{k-1}, l_k)$ be
a sequence of non-negative integers which is not $(0, 0, \dots, 0)$.
We call $L$ a $\lambda$-sequence if it satisfies
the following conditions.
\begin{enumerate}
\item $l_i\le q_i+1$ if $i\not= k$, and $l_k\le q_k$;
 $j_i\le q_i$ for all odd $i$, and $j_i\not= 1$ if $i\not= 1$
\item If $l_{i_0}=q_{i_0}+1$, then there exists odd $i_1$ and $i_2$
 which satisfies the following conditions.
 \begin{enumerate}
 \item $i_1<i_0<i_2\le k-1$
 \item $l_{i'}=q_{i'}$ for all even $i'$
 such that $i_1<i'<i_2$ and $i'\not= i_0$;
 $j_{i'}=0$ for all odd $i'$ such that $i_1\le i'\le i_2$
 \item $l_{i_1-1}<q_{i_1-1}$ if $i_1\ge 3$ and $l_{i_2+1}<q_{i_2}$
 \end{enumerate}
\item If $l_{i_0}=q_{i_0}$ and $j_{i_0+1}\ge 2$, then there exists odd $i_3$
  which satisfies the following conditions.
 \begin{enumerate}
 \item $i_3<i_0$
 \item $l_{i'}=q_{i'}$ for all even $i'$ such that $i_3<i'\le i_0$;
 $j_{i'}=0$ for all odd $i'$ such that $i_3\le i'<i_0$
 \item $l_{i_3-1}<q_{i_3-1}$ if $i_3\ge 3$
 \end{enumerate}
\end{enumerate}
\end{dfn}

\begin{dfn}
Let $M=(l_1, j_2, l_3, j_4, \dots, l_{k-1}, j_k)$ be
a sequence of non-negative integers which is not $(0, 0, \dots, 0)$.
We call $M$ a $\mu$-sequence if it satisfies
the following conditions.
\begin{enumerate}
\item $l_i\le q_i+1$ for all odd $i$;
  $2\le j_i\le q_i$ for all even $i$
\item If $l_{i_0}=q_{i_0}+1$,
 then there exists even $i_1$ and $i_2$ which satisfies
  the following conditions.
 \begin{enumerate}
 \item $0\le i_1<i_0<i_2\le k$
 \item $l_{i'}=q_{i'}$ for all odd $i'$
  such that $i_1<i'<i_2$ and $i'\not= i_0$;
   $j_{i'}=0$ for all even $i'$ such that $i_1\le i'\le i_2$
 \item $l_{i_1-1}<q_{i_1}$ if $i_1\ge 2$,
   $l_{i_2+1}<q_{i_2+1}$ if $i_2\le k-2$
 \end{enumerate}
\item If $l_{i_0}=q_{i_0}$ and $j_{i_0+1}\ge 2$, then there exists even $i_3$
  which satisfies the following conditions.
 \begin{enumerate}
 \item $0\le i_3<i_0$
 \item $l_{i'}=q_{i'}$ for all odd $i'$ such that $i_3<i'<i_0$;
   $j_{i'}=0$ for all even $i'$ such that $i_3\le i'<i_0$
 \item $l_{i_3-1}<q_{i_3-1}$ if $i_3\ge 2$
 \end{enumerate}
\end{enumerate}
\end{dfn}

We denote by $\cal S_{\lambda}$ (resp. $\cal S_{\mu}$)
the set of all $\lambda$- (resp. $\mu$-) sequences.\\
Let $L$ be a $\lambda$-sequence and $h$ an integer
such that $1\le h\le k-1$.
We say that {\it the condition $*(h)$ holds for $L$}
if the following condions hold.\\
\begin{enumerate}
\item If $h$ is odd,
$$\lambda^1_{j_1}+\sum
\begin{Sb}
3\le i\le h \\ i\ \text{odd}
\end{Sb}
(l_{i-1}\lambda^i_1+\lambda^i_{j_i})
<P_h+Q_h$$
\item If $h$ is even,
$$\lambda^1_{j_i}+\sum
\begin{Sb}
3\le i\le h-1 \\ i\ \text{odd}
\end{Sb}
(l_{i-1}\lambda^i_1+\lambda^i_{j_i})
+l_h\lambda^{h+1}_1
<P_{h-1}+Q_{h-1}+P_h+Q_h$$
If $l_h<q_h$ or $j_{h+1}\ge 2$ also hold,
$$\lambda^1_{j_1}+\sum
\begin{Sb}
3\le i\le h-1 \\ i\ \text{odd}
\end{Sb}
(l_{i-1}\lambda^i_1+\lambda^i_{j_i})
+l_h\lambda^{h+1}_1
<P_h+Q_h$$
\end{enumerate}

\begin{lem}
$*(h)$ holds for all $L\in \cal S _{\lambda}$
and for all $h=1, 2, \dots, k-1$.
\end{lem}

\begin{pf}
We use the induction on $h$.
It is clear that $*(1)$ holds.
Assume that $2\le h\le k-1$ and that
$*(\tilde{h})$ holds for all $\tilde{h}$
such that $\tilde{h} <h$.\\
First we treat that the case where $h$ is odd.\\
If $j_h\ge 2$, then by $*(h)$,
$$\lambda^1_{j_1}+\sum
\begin{Sb}
3\le i\le h-2 \\ i\ \text{odd}
\end{Sb}
(l_{i-1}\lambda^i_1+\lambda^i_{j_i})
l_{h-1}\lambda^h_1
<P_{h-1}+Q_{h-1}$$
Hence
\begin{align*}
\lambda^1_{j_1}+\sum
\begin{Sb}
3\le i\le h \\ i\ \text{odd}
\end{Sb}
(l_{i-1}\lambda^i_1+\lambda^i_{j_i})
&<P_{h-1}+Q_{h-1}+\lambda^h_{j_h}\\
&=P_{h-1}+Q_{h-1}+j_h(P_{h-1}+Q_{h-1})\\
&\le P_{h-2}+Q_{h-2}+q_h(P_{h-1}+Q_{h-1})\\
&=P_h+Q_h
\end{align*}

If $j_h=0$, by $*(h-1)$
\begin{align*}
\lambda^1_{j_1}+\sum
\begin{Sb}
3\le i\le h \\ i\ \text{odd}
\end{Sb}
(l_{i-1}\lambda^i_1+\lambda^i_{j_i})
&=\lambda^1_{j_1}+\sum
\begin{Sb}
3\le i\le h-2 \\ i\ \text{odd}
\end{Sb}
(l_{i-1}\lambda^i_1+\lambda^i_{j_i})+l_{h-1}\lambda^h_1\\
&<P_{h-2}+Q_{h-2}+P_{h-1}+Q_{h-1}\\
&\le P_h+Q_h
\end{align*}

Next we treat the case where $h$ is even.
We divide the proof into four cases as follows\\
(1) $l_h<q_h$ (2) $l_h=q_h, j_h=0$
(3) $l_h=q_h, j_h\ge 2$ (4) $l_h=q_h+1$\\
(1) By $*(h-1)$,
$$\lambda^1_{j_1}+\sum
\begin{Sb}
3\le i\le h-1 \\ i\ \text{odd}
\end{Sb}
(l_{i-1}\lambda^i_1+\lambda^i_{j_i})
<P_{h-1}+Q_{h-1}$$
Hence
\begin{align*}
\lambda^1_{j_1}+\sum
\begin{Sb}
3\le i\le h-1 \\ i\ \text{odd}
\end{Sb}
(l_{i-1}\lambda^i_1+\lambda^i_{j_i})
+l_k\lambda^{h+1}_1
&<P_{h-1}+Q_{h-1}+(q_h-1)(P_{h-1}+Q_{h-1})\\
&<P_h+Q_h
\end{align*}

(2) By $*(h-1)$,
\begin{align*}
\lambda^1_{j_1}+\sum
\begin{Sb}
3\le i\ h-1 \\ i\ \text{odd}
\end{Sb}
(l_{i-1}\lambda^i_1+\lambda^i_{j_i})
+l_h\lambda^{h+1}_1
&<P_{h-1}+Q_{h-1}+q_h(P_{h-1}+Q_{h-1})\\
&<P_{h-1}+Q_{h-1}+P_h+Q_h
\end{align*}

(3) By (iii) in Definition \ref{dfnlambda}, there exists odd $h'$
  such that $l_{h'-1}<q_{h'-1}$ and
$$(j_{h'}, l_{h'+1}, j_{h'+2}, l_{h'+3}, \dots, j_{h-1}, l_h)
=(0, q_{h'+1}, 0, q_{h'+3}, \dots, 0, q_h)$$
Hence by $*(h'-1)$,
\begin{align*}
&\lambda^1_{j_1}+\sum
\begin{Sb}
3\le i\le h-1 \\ i\ \text{odd}
\end{Sb}
(l_{i-1}\lambda^i_1+\lambda^i_{j_i})
l_h\lambda^{h+1}_1\\
=&\lambda^1_{j_1}+\sum
\begin{Sb}
3\le i\le h'-2 \\ i\ \text{odd}
\end{Sb}
(l_{i-1}\lambda^i_1+\lambda^i_{j_i})
+l_{h'-1}\lambda^{h'}_1\\
&+\lambda^{h'}_{j_{h'}}+\sum
\begin{Sb}
h'+2\le i\le h-1 \\ i\ \text{odd}
\end{Sb}
(l_{i-1}\lambda^i_1+\lambda^i_{j_i})
+l_{h}\lambda^{h+1}_1\\
<&P_{h'-1}+Q_{h'-1}+\sum
\begin{Sb}
h'+2\le i\le h+1 \\ i\ \text{odd}
\end{Sb}
q_{i-1}(P_{i-2}+Q_{i-2})\\
=&P_h+Q_h
\end{align*}

(4) By (ii) in Definition \ref{dfnlambda}, there exists odd $h'$ such that
  $l_{h'-1}<q_{h'-1}$ and
$$(j_{h'}, l_{h'+1}, j_{h'+2}, l_{h'+3}, j_{h'+4},
\dots, l_{h-2}, j_{h-1}, l_h)$$
$$=(0, q_{h'+1}, 0, q_{h'+3}, 0, \dots, q_{h-2}, 0, q_h)$$
Hence by $*(h'-1)$,
\begin{align*}
&\lambda^1_{j_1}+\sum
\begin{Sb}
3\le i\le h-1 \\ i\ \text{odd}
\end{Sb}
(l_{i-1}\lambda^i_1+\lambda^i_{j_i})
+l_h\lambda^{h+1}_1\\
=&\lambda^1_{j_1}+\sum
\begin{Sb}
3\le i\le h'-2 \\ i\ \text{odd}
\end{Sb}
(l_{i-1}\lambda^i_1+\lambda^i_{j_i})
+l_{h'-1}\lambda^{h'}_1\\
&+\lambda^{h'}_{j_{h'}}+\sum
\begin{Sb}
h'+2\le i\le h-1 \\ i\ \text{odd}
\end{Sb}
(l_{i-1}\lambda^i_1+\lambda^i_{j_i})
+l_h\lambda^{h+1}_1 \\
<&P_{h'-1}+Q_{h'-1}+\sum
\begin{Sb}
h'+2\le i\le h-1 \\ i\ \text{odd}
\end{Sb}
q_{i-2}(P_{i-2}+Q_{i-2})
(q_h+1)(P_{h-1}+Q_{h-1})\\
=&P_{h-1}+Q_{h-1}+P_h+Q_h
\end{align*}
\end{pf}

We define the order in $\cal S_{\lambda}$ as
$(j_1, l_2, \dots, l_k)<(j'_1, l'_2, \dots, l'_k)$
if and only if there exists $i$ such that
$j_i<j'_i$ or $l_i<l'_i$ and that $j_h=j'_h$, $l_h=l'_h$
for all $h>i$.

Let $v$ be a map from $\cal S_{\lambda}$ to $\Bbb Z$ defined by
$$v(l_1, j_2, \dots, l_{k-1}, j_k)
=\lambda^1_{j_1}+\sum
\begin{Sb}
3\le i\le k-1 \\ i\ \text{odd}
\end{Sb}
(l_{i-1}\lambda^i_1+\lambda^i_{j_i})
+l_k\lambda^{k+1}_1$$

\begin{prop}\label{proplambda}
The map $v$ is an order isomorphism from $\cal S_{\lambda}$
to $\{ n\in \Bbb Z |1\le n\le m-1\}$
\end{prop}

\begin{pf}
First we show that $L<L'$ implies $v(L)<v(L')$.
Let $L=(j_1, l_2, \dots, j_{k-1}, l_k)$ and
$L'=(j'_1, l'_2, \dots, j'_{k-1}, l'_k)$
be $\lambda$-sequences such that $L<L'$.
Put
$$i_0=\max\{ i|
l_{\tilde{i}}=l'_{\tilde{i}}\ and\ j_{\tilde{i}}=j'_{\tilde{i}}
\ \text{for all}\ \tilde{i} >i\}$$
If $i_0$ is odd,
$$v(L')-v(L)
\ge \lambda^{i_0}_{j'_{i_0}}-\lambda^{i_0}_{j_{i_0}}
-\{\lambda^1_{j_1}+\sum
\begin{Sb}
3\le i\le i_0-2 \\ i\ \text{odd}
\end{Sb}
(l_{i-1}\lambda^i_1+\lambda^i_{j_i})
+l_{i_0-1}\lambda^{i_0}_1\}$$
Note that
\begin{equation*}
\lambda^{i_0}_{j'_{i_0}}-\lambda^{i_0}_{j_{i_0}}\ge
\begin{cases}
P_{i_0-1}+Q_{i_0-1} &(\text{if}\ j_{i_0}\ge 2)\\
P_{i_0-2}+Q_{i_0-2}+P_{i_0-1}+Q_{i_0-1} &(\text{if}\ j_{i_0}=0)
\end{cases}
\end{equation*}
Hence by $*(i_0-1)$, $v(L')-v(L)>0$\\
If $i_0$ is even,
\begin{align*}
v(L')-v(L)
&\ge (l'_{i_0}-l_{i_0})\lambda^{i_0+1}_1
-\{ \lambda^1_{j_1}+\sum
\begin{Sb}
3\le i\le i_0-1 \\ i\ \text{odd}
\end{Sb}
(l_{i-1}\lambda^i_1+\lambda^i_{j_i})\}\\
&\ge P_{i_0-1}+Q_{i_0-1}-\{\lambda^1_{j_1}+\sum
\begin{Sb}
3\le i\le i_0-1 \\ i\ \text{odd}
\end{Sb}
(l_{i-1}\lambda^i_1+\lambda^i_{j_i})\}
\end{align*}

Hence by $*(i_0-1)$, $v(L')-v(L)>0$. Thus we have done.\\
Note that
$\max\cal S_{\lambda}
=(0, q_2, 0, q_4, 0, \dots, q_{k-2}, 0, q_k)$.
Thus $v$ is an order-preserving injection into
$\{ n\in\Bbb Z |1\le n\le m-1\}$.
Hence the rest we must prove is that
it is a injection into $\{ n\in\Bbb Z |1\le n\le m-1\}$.
Note that $1=\lambda^1_1<\lambda^1_2<\dots <\lambda^1_{q_1}
<\lambda^3_1<\lambda^3_2<\dots$.
Thus it is sufficient to show that
$\text{Im}(v)\supseteq \{n|1\le n<\lambda^h_1\}$
implies
$\text{Im}(v)\supseteq \{n|1\le n<\lambda^{h+2}_1\}$.
Let $n$ be an integer such that
$\lambda^h_1\le n<\lambda^{h+2}_1$.
We divide the proof into two cases.\\
(I) $\lambda^h_1\le n<\lambda^h_2$ or $q_h=1$
(II) $\lambda^h_2\le n<\lambda^{h+2}_1$\\
(I) Write $n=l_{h-1}\lambda^h_1+n'$ such that $0\le n'<\lambda^h_1$.
It is clear that $l_{h-1}\le q_{h-1}+1$.
By the induction hypothesis,
there exists $L'=(j_1, l_2, \dots, l_{h-3}, j_{h-2}, 0, \dots, 0)$
such that $v(L')=n'$.
Put $L=(j_1, l_2, \dots, l_{h-3}, j_{h-2}, l_{h-1}, 0, \dots, 0)$.
Assume that $L$ is not a $\lambda$-sequence.
Then by the definition of $\lambda$-sequence,
we get $l_{h-1}=q_{h-1}+1$ and there exists $h'$ such that
$$(l_{h'}, j_{h'+1}, l_{h'+2}, j_{h'+3}, l_{h'+4}, j_{h'+5},
\dots, l_{h-3}, j_{h-2}, l_{h-1})$$
$$=(q_{h'}+1, 0, q_{h'+2}, 0, q_{h'+4}, 0, \dots, q_{h-3}, 0, q_{h-1}+1)$$
or
$$(j_{h'}, l_{h'+1}, j_{h'+2}, l_{h'+3}, j_{h'+4},
\dots, l_{h-3}, j_{h-2}, l_{h-1})$$
$$=(j_{h'}, q_{h'+1}, 0, q_{h'+3}, 0, \dots, q_{h-3}, 0, q_{h-1}+1)$$
and $j_{h'}>0$.
In the both cases it is easily checked
\begin{equation*}
n\ge
\begin{cases}
\lambda^h_2 &(\text{if}\ q_h\ge 2)\\
\lambda^{h+2}_1 &(\text{if}\ q_h=1)
\end{cases}
\end{equation*}
and this is a contradiction.
Thus $L$ is a $\lambda$-sequence, so we have done.\\
(II) Put $j_h=\max\{ j|j\ge 2, \lambda^h_j\le n\}$.
Clearly $j_h\le q_h$.
Since $n-\lambda^h_{j_h}<q_h$,
there exists $L'=(j_1, l_2, \dots, j_{h-2}, l_{h-1}, 0, \dots, 0)$
such that $v(L')=n-\lambda^h_{j_h}$
by the induction hypothesis and (I).
Put $L=(j_1, l_2, \dots, j_{h-2}, l_{h-1}, j_h, 0, \dots, 0)$
We can check $L\in\cal S_{\lambda}$ by the definition of
$\lambda$-sequence.
\end{pf}

 When we write
$n=\lambda^1_{j_1}+\sum
\begin{Sb}
3\le i\le k-1 \\ i\ \text{odd}
\end{Sb}
(l_{i-1}\lambda^i_1+\lambda^i_{j_i})+l_k\lambda^{k+1}_1$
where $(j_1, l_2, \dots, l_{k-1}, j_k)$
is a $\lambda$-sequence,
we call this expression a {\it $\lambda$-expansion of $n$}.
By the above proposition, $n=1, 2, \dots, m-1$ has unique $\lambda$-
expansion.
Note that the proof of Proposition \ref{proplambda} shows
how to calculate $\lambda$-expansion of actual number.
\par
When we write
$n=\sum
\begin{Sb}
2\le i\le k\\i\ \text{even}
\end{Sb}
(l_{i-1}\mu^i_1+\mu^i_{j_i})$
where $(l_1, j_2, \dots, l_{k-1}, j_k)$,
we call this expression a {\it $\mu$-expansion of $n$}.
Similarly to $\lambda$-expansion,
we can prove that $n=1, 2, \dots, m-1$ has unique $\mu$-expansion.

\section{General $n$-canonical divisors}
\label{secdiv}
Let $(Y,y)$ be a 2-dimensional rational singularity, and $p:\yt\rightarrow Y$
be the minimal desingularization.
Let $M$ be a reflexive module of rank 1 on $Y$, $F(M)$ the full sheaf
associated to $M$. (For the definition of full sheaf, see \cite{esn}.
)
In this situation,
\begin{dfn}
Let $D$ be a member of $|M|$.
We call $D$ a general member of $|M|$ if \dt\ is a member of $|F(M)|$ and
intersects the exeptional locus transversely.
\end{dfn}
 Note that general members always exist since the full sheaf is generated by
global sections.
Let $E=\cup_i E_{i}$ be the exeptional locus,
and write $p^{*}D=\dt+\sum_i\alpha(D)_{i}E_{i}$.
\begin{lem}\label{lemfull}
 Let $D$ be a member of $|M|$ such that \dt\ and $E$ intersects transeversely.
 Then $D$ is a general member if and only if the inequality
$\alpha(D)_{i}\leq\alpha(D^{\prime})_{i}$ holds for all $D^{\prime}\in|M|$ and
all $E_{i}$.
\end{lem}
\begin{pf}
 Suppose that $D$ is a general member of $|M|$ and $D'$ is a member of $|M|$.
The sequence
$$H_{E}^{0}(\yt,\cal O_{\yt}(\dt))\rightarrow H^{0}(\yt,\cal
O_{\yt}(\dt))\rightarrow H^{0}(\yt\setminus E,\cal O_{\yt}(\dt))\rightarrow
H_{E}^{1}(\yt,\cal O_{\yt}(\dt))$$
is exact.
Since $\cal O_{\xt}(\dt)$ is a full sheaf, $H_{E}^{0}(\yt,\cal
O_{\yt}(\dt))=H_E^1(\yt,\cal O_{\yt}(\dt))=0$.
Hence
$H^0(\tilde{Y},\cal O_{\tilde{Y}}(\tilde{D}))\simeq H^0(Y, \cal O_Y(D))$.
Thus $D$ is linearly equivalent to $D'+(\text{effective divisor})$.
Hence the inequality holds since $\dt$ and
$\dt^{\prime}+\sum\{\alpha(D^{\prime})_i-\alpha(D)_i\}E_i$
are linearly equivalent. Thus we have proved only if part.
 Next suppose that $D$ is a member of $|M|$ such that \dt\ and $E$ intersect
transeversely and $\alpha(D)_i\leq\alpha(D^{\prime})_i$ for all
$D^{\prime}\in|M|$.
Let $D_0$ be a general member of $|M|$. By the assumption of $D$,
$\alpha(D)_i\leq\alpha(D_0)_i$.
By the fact that we have already showed, $\alpha(D)_i\geq\alpha(D_0)_i$. Hence
$\alpha(D)_i=\alpha(D_0)_i$.
 Thus $D$ and $D_0$ are numerically equivalent,
hence they are lenearly equivalent since $Y$ is a rational singularity.
Hence $\cal O_{\yt}(\dt)$ is a full sheaf.
\end{pf}
\begin{cor}\label{corgenint}
Let $M$ and $M^{\prime}$ be reflexive modules of rank 1 on $Y$. Let
$D_0$(resp.$D_0^{\prime}$) be a general member of $|M|$(resp.$|M^{\prime}|$).
Then
$$D_0\cdot D_0^{\prime}=\min\{D\cdot
D^{\prime}|D\in|M|,D^{\prime}\in|M^{\prime}|\}$$
\end{cor}
\begin{pf}
Let $D$(resp.$D^{\prime}$) be a member of $|M|$(resp.$|M^{\prime}|$).
\begin{align*}
D\cdot D^{\prime}  &
=(p^*D)\cdot\dt^{\prime} \\
 & =(\dt+\sum\alpha(D)_iE_i)\cdot\dt^{\prime} \\
 & \geq\sum\alpha(D)_i(E_i\cdot\dt^{\prime}) \\
 & \geq\sum\alpha(D_0)_i(E_i\cdot\dt^{\prime}) \\
 & =(\dt_0+\sum\alpha(D_0)_iE_i)\cdot\dt^{\prime} \\
 & =p^*D_0\cdot\dt_0=D_0\cdot D^{\prime}
\end{align*}
Similarly we can show $D_0\cdot D'\ge D_0\cdot D_0'$.
\end{pf}

 In the rest of this section we shall calculate
the general element of the $n$-canonical system of
semi-log-terminal singularities.
We denote by $\cal L$ (resp. $\cal L_o$ resp. $\cal L_e$)
the set of all functions from $I$ (resp. $I_o$ resp. $I_e$)
to $\Bbb Z$.
\par
First we treat singularities of class T.
We begin by purely arithmetical lemmas.
 Let $a$, $d$, $m$ and $n$ be positive integers
such that $\dfrac{m}{2}<a<m$, $(a, m)=1$ and $n<m$.
Let $\dfrac{a}{m-a} =[q_1, q_2, \dots, q_k]$
be the expantion into continued fraction,
$r_i$ be the $i$-th remainder of the Euclidean algorithm,
and $\dfrac{P_i}{Q_i}$ be the $i$-th convergent.

\begin{lem}\label{lemrineq}
Let $\{ t_i\}_{i=1, 2, \dots, k}$ be a sequence of non-negative integers
such that $t_i\le q_i$ for $i\le k-1$ and $t_k\le q_k-1$.
Assume that $t_{i_0}$ be positive. Then
\begin{enumerate}
\item If $i_0$ is even,
$$-r_{i_0-1}<\sum_{i_0\le i\le k}(-1)^{i-1}t_ir_i<0$$
\item If $i_0$ is odd,
$$0<\sum_{i_0\le i\le k}(-1)^{i-1}t_ir_i<r_{i_0-1}$$
\end{enumerate}
\end{lem}

\begin{pf}
We use the induction on $k-i_0$.
If $i_0=k$, it is clear that the inequalities hold.
Assume that $i_0<k$ and that the inequalities hold
for all $i'_0$ such that $i_0<i'_0\le k$.
Let $i_0$ be even.
First we show $\sum_{i_0\le i\le k}(-1)^{i-1}t_ir_i<0$.
If $t_{i'}=0$ for all $i'$ such that $i'>i_0$ and $i'$ is odd,
$$\sum_{i_0\le i\le k}(-1)^{i-1}t_ir_i\le -t_{i_0}r_{i_0}<0$$
Thus we assume that there exists $i'$ such that $i'>i_0$, $i'$ is odd, and
$t_{i'}>0$.
Let $i'_0$ be the minimum of such $i'$.
By the induction hypothesis,
$$\sum_{i'_0\le i\le k}(-1)^{i-1}t_ir_i<r_{i'_0}-1$$
Thus
\begin{align*}
\sum_{i'_0\le i\le k}(-1)^{i-1}t_ir_i
&\le -t_{i_0}r_{i_0}+\sum_{i'_0\le i\le k}(-1)^{i-1}t_ir_i\\
&<-t_{i_0}r_{i_0}+r_{i'_0-1}\\
&=0
\end{align*}

Next we show $\sum_{i_0\le i\le k}(-1)^{i-1}t_ir_i>-r_{i_0-1}$.
If $t_{i'}=0$ for all $i'$ such that $i'>i_0$ and $i'$ is even,
$$\sum_{i_0\le i\le k}(-1)^{i-1}t_ir_i\ge -t_{i_0}r_{i_0}
\ge -q_{i_0}r_{i_0}>-r_{i_0-1}$$
Thus we assume that there exists $i'$ such that $i'>i_0$, $i'$ is even, and
$t_{i'}>0$.
By the induction hypothesis,
$$\sum_{i'_0\le i\le k}(-1)^{i-1}t_ir_i>-r_{i'_0-1}$$
Thus
\begin{align*}
\sum_{i_0\le i\le k}(-1)^{i-1}t_ir_i
&\ge -t_{i_0}+\sum_{i_0\le i\le k}(-1)^{i-1}t_ir_i\\
&>-t_{i_0}-r_{i'_0-1}\\
&\ge -q_{i_0}r_{i_0}-r_{i_0+1}\\
&=-r_{i_0-1}
\end{align*}

The proof is similar for odd $i_0$.
\end{pf}

\begin{dfn}
Let $(t_1, t_2, \dots, t_k)$ be a sequence of non-negative integers which
is not $(0, 0, \dots, 0)$.
We call this a $\tau$-sequence if it satisfies the following conditions.
\begin{enumerate}
\item $t_i\le q_i$ if $i\not= k$ and $t_k<q_k$.
\item If $t_{i_0-1}>0$ and $t_{i_0}=q_{i_0}$
 for some $i_0$ such that $1<i_0<k$, then $t_{i_0+1}=q_{i_0+1}$.
    If $t_{k-2}>0$ and $t_{k-1}=q_{k-1}$, then $t_k=q_k-1$.
\end{enumerate}
\end{dfn}

\begin{lem}\label{lemtau}
Let $(0, 0, \dots, 0, t_{i_0}, t_{i_0+1}, \dots, t_k)$ be a $\tau$-sequence
such that $i_0$ is odd and $t_{i_0}>0$. Then
\begin{enumerate}
\item If $k$ is even,
$$\sum_{i_0\le i\le k}(-1)^{i-1}t_ir_i=1\Leftrightarrow i_0=k-1, t_{k-1}=1,
t_k=q_k-1$$
\item If $k$ is odd,
$$\sum_{i_0\le i\le k}(-1)^{i-1}t_ir_i=1\Leftrightarrow i_0=k, t_k=1$$
\end{enumerate}
\end{lem}

\begin{pf}
It can be easily checked and we left it for the reader.
\end{pf}

\begin{lem}
For any integer $t$ such that $0<t<m-2$,
there exists a $\tau$-sequence $(t_1, t_2, \dots, t_k)$ such that
$t=\sum_{1\le i\le k}t_i(P_{i-1}+Q_{i-1})$.\\
(We call this expression of $t$ the $\tau$-expansion of $t$.)
\end{lem}

\begin{pf}
(Step 1) We can write $t=\sum_{1\le i\le k}t_i(P_{i-1}+Q_{i-1})$
where $0\le t_i\le q_i$ for $i<k$ and $0\le t_k<q_k$.\\
(proof) We use the induction on $t$. By the induction hypothesis,
we can write $t-1=\sum_{1\le i\le k}t_i(P_{i-1}+Q_{i-1})$,
where $0\le t_i\le q_i$ for $i<k$ and $0\le t_k<q_k$.
If $t_1<q_1$, we have done.
Assume $t_1=q_1$.
Note that $\sum_{1\le i\le
k-1}q_i(P_{i-1}+Q_{i-1})+(q_k-1)(P_{k-1}+Q_{k-1})=m-2$.
Hence there exists $i_0$ such that $i_0<k$, $t_i=q_i$ for all $i\le i_0$,
and $t_{i_0+1}<q_{i_0+1}$ (if $i_0<k-1$); $t_k\le q_k-2$ (if $i_0=k-1$).
Thus
\begin{equation*}
t=
\begin{cases}
\sum
\begin{Sb}
1\le i\le i_0 \\ i\ \text{even}
\end{Sb}
q_i(P_{i-1}+Q_{i-1})+(t_{i_0+1}+1)(P_{i_0}+Q_{i_0})
+&\sum_{i_0+2\le i\le k}t_i(P_{i-1}+Q_{i-1}) \\
&(\text{if}\ i_0\ \text{is odd})\\
\sum
\begin{Sb}
1\le i\le i_0 \\ i\ \text{odd}
\end{Sb}
q_i(P_{i-1}+Q_{i-1})+(t_{i_0+1}+1)(P_{i_0}+Q_{i_0})
+&\sum_{i_0+2\le i\le k}t_i(P_{i-1}+Q_{i-1}) \\
&(\text{if}\ i_0\ \text{is even})
\end{cases}
\end{equation*}

(Step 2) We write $t=\sum_{1\le i\le k}t^{(1)}_i(P_{i-1}+Q_{i-1})$ as in (Step
1).
If there exists $i_0$ such that $t_{i_0-1}>0$, $t_{i_0}=q_{i_0}$,
and $t_{i_0+1}<q_{i_0+1}$ (if $i_0<k-1$); $t_{i_0}\le q_{i_0+1}-2$ (if
$i_0=k-1$),
we transform $(t^{(1)}_1, \dots, t^{(1)}_k)$ to
$(t^{(2)}_1, \dots, t^{(2)}_k)
=(t_1, \dots, t_{i_0-2}, t_{i_0-1}-1, 0, t_{i_0+1}+1, t_{i_0+2}, \dots, t_k)$.
We repeat this operation.
Since $\sum_{1\le i\le k}t^{(j)}_i$ strictly decrease by this operation,
we can get a $\tau$-sequence after finitely many operations.
\end{pf}

Remark. We can prove that the $\tau$-expansion is unique.

Put $\cal T$, $v$ and $\cal T_{\min}$ as follows
\begin{align*}
\cal T &=\{ (s, t)\in \Bbb Z _{\ge 0} \times \Bbb Z _{\ge 0}|
 s+\{ dm(m-a)-1\} t\equiv dm(m-a)n\ (mod\ dm^2) \}\\
v&=\min\{s+t|(s, t)\in \cal T \}\\
\cal T _{\min}&=\{ (s, t)\in \cal T|s+t=v\}
\end{align*}

\begin{prop}\label{propsmin}
\begin{enumerate}
\item If $d=1$ and $k$ is even, then
\begin{equation*}
\begin{cases}
 \cal T _{\min}=\{ (n, n)\} &(n<m-(P_{k-1}+Q_{k-1}))\\
 \cal T _{\min}=\{ (n+P_{k-1}+Q_{k-1}, n-m+P_{k-1}+Q_{k-1})\}
 \ &(n\ge m-(P_{k-1}+Q_{k-1}))
\end{cases}
\end{equation*}

\item If $d=1$ and $k$ is odd, then
\begin{equation*}
\begin{cases}
 \cal T _{\min}=\{ (n, n)\} &(n<m-(P_{k-1}+Q_{k-1}))\\
 \cal T _{\min}=\{ (n-m+P_{k-1}+Q_{k-1}, n+P_{k-1}+Q_{k-1})\}
 \ &(n\ge m-(P_{k-1}+Q_{k-1}))
\end{cases}
\end{equation*}

\item If $d\ge 2$, then $\cal T_{\min} =\{ (n, n)\}$
\end{enumerate}
\end{prop}

\begin{pf}
For a positive integer $t$, put
$$s_t=\min\{ s|(s, t)\in \cal T \}$$
Since $(n, n)\in \cal T$, we get $v\le 2n$, thus
$$v=\min\{ s_t+t|0\le t\le 2n\}$$
$$\cal T _{\min}=\{(s_t, t)|0\le t\le 2n, s_t+t=v\}$$

Hence we estimate $s_t$ for $t$ such that $0\le t\le 2n$.
It is clear that $s_n=n$.\\
(Claim I) Assume $t<n$. Then
\begin{enumerate}
\item If $k$ is even, $d=1$, $n\ge m-(P_{k-1}+Q_{k-1})$
 and $t=n-m+(P_{k-1}+Q_{k-1})$,
 then $s_t=n+P_{k-1}+Q_{k-1}$
\item Otherwise $s_t+t>2n$
\end{enumerate}
(Proof of Claim I)
It is clear that $s_0=dma>2n$, thus we assume $t>0$.
Put $t'=n-t$.
Let $t'=\sum_{1\le i\le k}t_i(P_{i-1}+Q_{i-1})$
be the $\tau$-expansion of $t'$.
Put $i_0=\min\{i|t_i>0\}$
First let $i_0$ be even.
We must prove $s_t+t>0$ in this case
since $i_0=1$ when $k$ is even and $t=n-m+P_{k-1}+Q_{k-1}$.
Let $(s, t)$ be an element in $\cal L$.
We can write $s=n+\{dm(m-a)-1\}t'-s'$ in which $s'$ is an integer.
Note that $(m-a)t'-m\sum_{1\le i\le k}t_iQ_{i-1}
=\sum_{1\le i\le k}(-1)^{i-1}t_ir_i$.
Assume $s'\ge \sum_{1\le i\le k}t_iQ_i$.
Then by Lemma \ref{lemrineq},

\begin{align*}
s&\le n-t'+dm\{ (m-a)t'-m\sum_{1\le i\le k}t_iQ_{i-1}\}\\
 &=n-t'+dm\sum_{1\le i\le k}(-1)^{i-1}t_ir_i\\
 &\le n-t'-dm\\
 &<0
\end{align*}

This is a contradiction.
Thus $s'\le \sum_{1\le i\le k}t_iQ_{i-1}-1$.\par
Next put $s'=\sum_{1\le i\le k}t_iQ_{i-1}-1$.
Then by Lemma \ref{lemrineq},
\begin{align*}
s&=n-t'+dm\sum_{1\le i\le k}(-1)^{i-1}t_ir_i+dm^2\\
 &>n-t'+dma\\
 &>2n
\end{align*}

Hence $s'_t=\sum_{1\le i\le k}t_iQ_{i-1}+1$ and $s_t>2n$.
Thus we have done for even $i_0$.\\

Secondly let $i_0$ be odd.
Assume that $s'\ge \sum_{1\le i\le k}t_iQ_{i-1}+1$.
Then by Lemma \ref{lemrineq},

\begin{align*}
s&\le n-t'+dm\sum_{1\le i\le k}(-1)^{i-1}t_ir_i-dm^2\\
 &<n-t'+dma-dm^2\\
 &<0
\end{align*}

This is a contradiction,
thus $s'\le \sum_{1\le i\le k}t_iQ_{i-1}+1$.
Next assume
$s'=\sum_{1\le i\le k}t_iQ_{i-1}$.
Then by Lemma \ref{lemrineq},
$$s=n-t'+dm\sum_{1\le i\le k}(-1)^{i-1}t_ir_i\ge n-t'+dm>0$$
Hence $$s_t=n-t'+dm\sum_{1\le i\le k}t_ir_i$$
and $$s_t+t=2t+dm\sum_{1\le i\le k}(-1)^{i-1}t_ir_i$$
If $d\sum_{1\le i\le k}(-1)^{i-1}t_ir_i\ge 2$, then $s_t+t=2m>2n$.
If $d\sum_{1\le i\le k}(-1)^{i-1}t_ir_i=1$, by Lemma \ref{lemtau},\\
(a) $d=1$, $k$ is even, and $t'=m-(P_{k-1}+Q_{k-1})$.\\
or (b) $d=1$, $k$ is odd, and $t'=P_{k-1}+Q_{k-1}$.\\
In the case (b),
$$s_t+t=2n-2(P_{k-1}+Q_{k-1})+m>2n$$
In the case (a),
$$t=n-m+P_{k-1}+Q_{k-1},\ s_t=n+P_{k-1}+Q_{k-1}$$
Thus we have done.

(Claim II) Assume $t>n$. Then
\begin{enumerate}
\item If $k$ is odd, $d=1$, $n\ge m-(P_{k-1}+Q_{k-1})$,
and $t=n+P_{k-1}+Q_{k-1}$,
then $s_t=n-m+P_{k-1}+Q_{k-1}$
\item Otherwise $s_t+t>2n$
\end{enumerate}
(Proof of Claim II)
Put $t'=t-n$.
If $t'=m-1$, then $n=m-1$ and $t=2m-2$.
It is easy to see
$$s_{2m-2}=dm(m-a)+m-1>2n$$
Thus we assume $t'\le m-2$.
Let $t'=\sum_{1\le i\le k}t_i(P_{i-1}+Q_{i-1})$
be the $\tau$-expansion of $t'$,
and put $i_0=\min\{i|t_i>0\}$.
Write $s=n-\{ dm(m-a)-1\}t'+dm^2s'$ for $(s, t)\in \cal T$.

For the case such that $i_0$ is even, we can get
$$s_t=n+t'-dm\sum_{1\le i\le k}(-1)^{i-1}t_ir_i>2n$$
by the similar way using Lemma \ref{lemrineq}.
If $t=P_{k-1}+Q_{k-1}$, then $i_0=k$, thus we have done for this case.

Assume that $i_0$ is odd, and let $(s, t)\in \cal T$.
If $s'\le \sum_{1\le i\le k}t_iQ_{i-1}-1$,
we can get $s<0$ by Lemma \ref{lemrineq}.
Put $s'=\sum_{1\le i\le k}t_iQ_{i-1}$.
Then we get $s=n+t'-dm\sum_{1\le i\le k}(-1)^{i-1}t_ir_i$
If $d\sum_{1\le i\le k}(-1)^{i-1}t_ir_i\ge 2$,
then $s_t=n+t'-dm\sum_{1\le i\le k}(-1)^{i-1}t_ir_i+dm^2$,
hence $s_t+t>2n$.\\
If $d\sum_{1\le i\le k}(-1)^{i-1}t_ir_i=1$, by Lemma \ref{lemrineq}\\
(a) $d=1$, $k$ is even, $t'=m-(P_{k-1}+Q_{k-1})$\\
or (b) $d=1$, $k$ is odd, $t'=P_{k-1}+Q_{k-1}$\\
In the case (a), $n+t'-dm\sum_{1\le i\le k}(-1)^{i-1}t_ir_i
=n-(P_{k-1}+Q_{k-1})$\\
If $n<P_{k-1}+Q_{k-1}$,
then $s'_t\ge \sum_{1\le i\le k}t_iQ_{i-1}+1$,
thus $s_t+t>2n$.\\
If $n\ge P_{k-1}+Q_{k-1}$,
then $s_t=n-(P_{k-1}+Q_{k-1})$,\\
thus $s_t+t=2n+m-2(P_{k-1}+Q_{k-1})>2n$.\\
In the case (b),
$$n+t'-dm\sum_{1\le i\le k}(-1)^{i-1}t_ir_i
=n+P_{k-1}+Q_{k-1}-m$$
If $n<m-(P_{k-1}+Q_{k-1})$,
then $s_t=n+t'-dm\sum_{1\le i\le k}(-1)^{i-1}t_ir_i+dm^2$,\\
thus $s_t+t>2n$.\\
If $n\ge m-(P_{k-1}+Q_{k-1})$,
then $s_t=n-m+P_{k-1}+Q_{k-1}$.\\
Thus we have done.

 Summarizing (I) and (II), we have proved the proposition.
\end{pf}

\par
We define $\delta^{\iota}\in \cal L$ by $\delta^{\iota}_{\eta}=0$
for $\eta\not= \iota$ and $\delta^{\iota}_{\iota}=1$.
And for $\nu\in\cal L$, we define $\alpha(\nu)\in\cal L\otimes\Bbb Q$
by $\alpha(\nu)_{\eta}=\sum_{\iota\in I}\alpha^{\iota}_{\eta}\nu_{\iota}$.

\begin{dfn}
Let $n$ be an integer such that $1\le n\le m-1$.
Let
\begin{equation*}
n=\lambda^1_{j_1}+\sum
\begin{Sb}
3\le i\le k-1 \\ i\ \text{\em{odd}}
\end{Sb}
(l_{i-1}\lambda^i_1+\lambda^i_{j_i})+l_k\lambda^{k+1}_1
=\sum
\begin{Sb}
2\le i\le k\\i\ \text{\em{even}}
\end{Sb}
(l_{i-1}\mu^i_1+\mu^i_{j_i})
\end{equation*}
be the $\lambda$- and $\mu$- expansion of $n$.
We define $\nu(n)^o\in I_o$, $\nu(n)^e\in I_e$,
$\nu(n)\in I$ as follows:
\begin{align*}
\nu(n)^o=&\delta^{1, j_1}+\sum
\begin{Sb}
3\le i\le k-1\\i\ \text{\em{odd}}
\end{Sb}
(l_{i-1}\delta^{i, 1}+\delta^{i, j_i})+l_k\delta^{k+1, 1}\\
\nu(n)^e=&\sum
\begin{Sb}
2\le i\le k\\i\ \text{\em{even}}
\end{Sb}
(l_{i-1}\delta^{i, 1}+\delta^{i, j_i})\\
\nu(n)=&\nu(n)^o+\nu(n)^e
\end{align*}
\end{dfn}

\begin{thm}\label{thmdegfull1}
Let $n$ be an integer such that $1\le n\le m-1$. Then
$$\deg_{E_{\iota}}F(-nK_X)=\nu (n)_{\iota}\quad
 for\ all\ \iota\in I$$
\end{thm}

\begin{pf}
Let $\nu'(n)$ be an element of $\cal L$ such that
$\nu'(n)_{\iota}=\deg_{E_{\iota}}F(-K_X)$.
Put $d'=\lceil \dfrac{d}{2} \rceil$.
Put $I'_o$ and $I'_e$ as follows
\begin{align*}
&I'_o=\{ (i, j)\in I|i \text{ is odd and } j\le d'\text{ if }i=k+1\}\\
&I'_e=I\smallsetminus I'_o
\end{align*}

For $\nu\in \cal L$, define $s(\nu)$ and $t(\nu)$ as follows
$$s(\nu)=\sum_{\iota\in I'_o}\nu_{\iota}\lambda_{\iota},\quad
 t(\nu)=\sum_{\iota\in I'_e}\nu_{\iota}\mu_{\iota}$$

Let $\cal L (n)$ be the set of elements of $\cal L$
satisfying the following conditions\\
(i) $\nu_{\iota}\ge 0$  for all $\iota\in I$\\
(ii) $s(\nu )+\{ dm(m-a)-1\} t(\nu) \equiv dm(m-a)n \pmod{dm^2}$\\

By Lemma \ref{lemfull}, it is clear
that $\nu'(n)$ is an element of $\cal L (n)$ which is
characterized by the inequalities $\alpha (\nu'(n))_{\eta}
=\alpha (\nu )_{\eta}$ for all $\nu\in\cal L$ and
for all $\eta\in I$.
First we show
$$s(\nu'(n))=s(\nu (n)),\quad t(\nu'(n))=t(\nu (n))$$
For $\nu\in\cal L$,
\begin{align*}
dm^2\alpha (\nu)^{k+1}_{d'}
=&\mu^{k+1}_{d'}\sum_{\iota\in I'_o}\nu_{\iota}\lambda_{\iota}
 +\lambda^{k+1}_{d'}\sum_{\iota\in I'_e}\nu_{\iota}\mu_{\iota}\\
=&\{ (d-d'+1)m-(P_{k-1}+Q_{k-1})\}(s(\nu )+t(\nu ))\\
&-\{(d-2d'+2)m-2(P_{k-1}+Q_{k-1})\}t(\nu )\\
=&:\alpha (s(\nu)+t(\nu))-\beta t(\nu)
\end{align*}

If $d\ge 2$,
\begin{align*}
dm^2\alpha (\nu )^{k+1}_{d'+1}
=&\{ (d-d')m-(P_{k-1}+Q_{k-1})\}(s(\nu)+t(\nu))\\
 &+\{ (2d'-d)m+2(P_{k-1}+Q_{k-1})\} t(\nu)\\
=&:\gamma (s(\nu )+t(\nu ))+\delta t(\nu )
\end{align*}

Note that $\alpha$, $\beta$, $\gamma$, and $\delta$ are
all positive.
Thus by the Proposition \ref{propsmin} we have done for this case.\\
If $d=1$,
$$dm^2\alpha (\nu)^k_{q_k -1}=
\{m-2(P_{k-1}+Q_{k-1})\}(s(\nu)+t(\nu))
+4(P_{k-1}+Q_{k-1})t(\nu)$$

Thus we can use the same arguement as above.\\

Next we show $\nu (n)_{\eta}=\nu'(n)_{\eta}$ for $\eta\in I_o$
by induction.
Let $\eta$ be an element of $I_o$ which is not (1, 1).
Assume $\nu (n)_{\iota}=\nu'(n)_{\iota}$ for all $\iota\in I_o$
such that $\iota >\eta$.
For $\nu\in\cal L$,
\begin{align*}
dm^2\alpha (\nu )_{\eta^l}=&-dm^2\nu_{\eta}+(s(\nu)-
 \sum_{\iota >\eta}\nu_{\iota}\lambda_{\iota})\mu_{\eta^l}\\
&+(t(\nu)+\sum_{\iota >\eta}\nu_{\iota}\mu_{\iota})\lambda_{\eta^l}
\end{align*}

Thus $$\nu'(n)_{\eta}\ge \nu(n)_{\eta}$$
Note that $$\sum_{\iota\in I'_o, \iota\le\eta}\nu'(n)_{\iota}\lambda_{\iota}
=\sum_{\iota\in I'_o, \iota\le\eta}\nu (n)_{\iota}\lambda_{\iota}$$
Thus by the property of the $\lambda$-expansion,
we get $\nu'(n)_{\eta}=\nu (n)_{\eta}$.\\
The same arguement shows $\nu (n)_{\eta}=\nu' (n)_{\eta}$
for all $\eta\in I_e$.
\end{pf}

\begin{cor}
Let $n$ be an integer such that $1\le n\le m-1$. Let
$$n=\lambda^1_{j_1}+\sum
\begin{Sb}
3\le i\le k-1 \\ i\ \text{\em{odd}}
\end{Sb}
(l_{i-1}\lambda^i_1+\lambda^i_{j_i})+l_k\lambda^{k+1}_1
=\sum
\begin{Sb}
2\le i\le k \\ i\ \text{\em{even}}
\end{Sb}
(l_{i-1}\mu^i_1+\mu^i_{j_i})$$
be the $\lambda$- and $\mu$- expansion.\\
Then
$$C^1_{j_1}+\sum
\begin{Sb}
3\le i\le k-1 \\ i\ \text{\em{odd}}
\end{Sb}
(\sum_{1\le h_{i-1}\le l_{i-1}}C^i_{1, h_{i-1}}
+C^i_{j_i})
+\sum_{1\le h_k\le l_k}C^{k+1}_{1, h_k}$$
$$+\sum
\begin{Sb}
2\le i\le k \\ i\ \text{\em{even}}
\end{Sb}
(\sum_{1\le h_{i-1}\le l_{i-1}}C^i_{1, h_{i-1}}
+C^i_{j_i})$$
is a general member of $|-nK_X|$.

\end{cor}

 Next we treat the non-normal case. It is easier than the
normal case.
 \begin{thm}
Let $n$ be an integer such that $1\le n\le m-1$. Then
$$\deg_{E_{\iota}}F(-n(K_{X_o}+\Delta_o))
=\nu_o(n)_{\iota}\quad for\ all\ \iota\in I_o$$
$$\deg_{E_{\iota}}F(-n(K_{X_e}+\Delta_e))
=\nu_e(n)_{\iota}\quad for\ all\ \iota\in I_e$$
\end{thm}

\begin{pf}
We only show the first equality since the proof is
similar for the second one.
Let $\nu'_o(n)$ be an element of $\cal L_o$ such that
$\nu'_o(n)_{\iota}=\deg_{E_{\iota}}F(-n(K_{X_o}+\triangle_o))$
for $\iota\in I_o$.
Put $\sigma (\nu)=\sum_{\iota\in I_o}\nu_{\iota}\lambda_{\iota}$
for $\nu\in\cal L_o$, and put
$$\cal L_o(n)=\{\nu\in\cal L_o|\nu\ \text{is nef},
\sigma (\nu)\equiv n\pmod{m}\}$$
By Lemma \ref{lemfull}, $\nu'_o(n)$ is an element of
$\cal L_o(n)$ which is characterized by
$\alpha(\nu'_o(n))_{\eta}\le\alpha(\nu)_{\eta}$
for all $\nu\in\cal L_o$ and for all $\eta\in I_o$.
Since $m\alpha(\nu)^{k+1}_1=\sigma(\nu)$
and $\sigma(\nu_o(n))=n$,
we get $\sigma(\nu'_o(n))=\sigma(\nu_o(n))$.
Thus we can prove the theorem by the induction using the formula
$$m\alpha(\nu)_{\eta^l}=-m\nu_{\eta}+
\mu_{\eta^l}(\sigma(\nu)-
\sum_{\iota\in I_o, \iota >\eta}\nu_{\iota}\lambda_{\iota})
+\lambda_{\eta^l}\sum_{\iota\in I_o, \iota >\eta}
\nu_{\iota}\lambda_{\iota}$$
similarly to the proof of Theorem \ref{thmdegfull1}.
\end{pf}

\section{Local intersection number}

As the application of the result of Section \ref{secdiv}, we shall prove
the following Thorem \ref{thmbound} in this section.
\begin{dfn}
 For a sequence of positive integers $(L_1, L_2, \dots, L_J)$ and
a positive integer $N$, we define the sequence $(N_{-1}, N_0, N_1, \dots,
N_J)$ by $$N_{-1}=N_0=N,\ N_j=L_jN_{j-1}+N_{j-2}\ (1\le j\le J)$$
We define $B((L_1, L_2, \dots, L_J), N)$ by
$B((L_1, L_2, \dots, L_J), N)=N_J$.
\par For a pair of positive integers $(M, N)$, we define $B(M, N)$ by
$$B(M, N)=\max\{B((L_1, L_2, \dots, L_J), N)|
\sum_{1\le j\le J}L_j=M\}$$
\end{dfn}

\begin{thm}\label{thmbound}
Let $(X, x)$ be a 2-dimensional smoothable semi-log-terminal
singularity, and $n$ a positive integer.
Let $D$ and $D'$ be members in $|nK_X|$ which do not have common components.
Then
$$\text{\em{index}}(X, x)\le B(D\cdot D'+1, n)$$
\end{thm}

We define $\Bbb Z$-valued symmetric bilinear forms
$\cal O$ and $\cal E$ on $\cal L$ by
\begin{equation*}
\cal O(\delta^{\iota},\delta^{\eta})=
\begin{cases}
\lambda_{\iota}\bar{\lambda}_{\eta}\quad &(\iota,\eta\in I_o,\iota <\eta)\\
\bar{\lambda}_{\iota}\lambda_{\eta}\quad &(\iota,\eta\in I_o,\iota\ge\eta)\\
0\quad &(\text{otherwise})
\end{cases}
\end{equation*}
\begin{equation*}
\cal E(\delta^{\iota},\delta^{\eta})=
\begin{cases}
\mu_{\iota}\bar{\mu}_{\eta}\quad &(\iota,\eta\in \bar{I}_e,\iota<\eta)\\
\bar{\mu}_{\iota}\mu_{\eta}\quad &(\iota,\eta\in \bar{I}_e,\iota\ge\eta)\\
0\quad &(\text{otherwise})
\end{cases}
\end{equation*}
where $\bar{I}_e=I_e\cup\{(k+1,d)\}$ in the normal case, and
$\bar{I}_e=I_e$ in the non-normal case.

\begin{lem}\label{lemint}
Let $\nu=\nu^o+\nu^e$ and $\tilde{\nu}=\tilde{\nu}^o+\tilde{\nu}^e$ be members
in $\cal L$.
Then in the normal case,
\begin{align*}
dm^2(\nu\cdot\tilde{\nu})
=&(dma-1)\sigma (\nu^o)\sigma (\tilde{\nu}^o) +\sigma (\nu^o)\tau
(\tilde{\nu}^e)
+\tau (\nu^e)\sigma (\tilde{\nu}^o) -(dma+1)\tau (\nu^e) \tau(\tilde{\nu}^e)\\
&+dm^2 (\cal E (\nu^e , \tilde{\nu}^e) -\cal O (\nu^o, \tilde{\nu}^o ))
\end{align*}
and in the non-normal case,
$$
m(\nu\cdot\tilde{\nu})
=a(\sigma(\nu^o)\sigma(\tilde{\nu}^o)
-\tau(\nu^e)\tau(\tilde{\nu}^e))
+m(\cal E(\nu^e,\tilde{\nu}^e)
-\cal O(\nu^o,\tilde{\nu}^o))
$$
\end{lem}

\begin{pf}
First, we treat the normal case.
We shall calculate the each term of the left hand side of
$$\nu \cdot \tilde{\nu} = \nu^o \cdot \tilde{\nu}^o +\nu^o \cdot \tilde{\nu}^e
+ \nu^e \cdot \tilde{\nu}^o + \nu^e \cdot \tilde{\nu}^e$$
We can easily get $dm^2 \nu^o \cdot \tilde{\nu}^e =\sigma (\nu^o)\tau
(\tilde{\nu}^e)$ and
$dm^2 \nu^e \cdot \tilde{\nu}^o = \tau (\nu^e)\sigma(\tilde{\nu}^o)$.
By the formula $\mu_{\iota}=-\lambda_{\iota}+dm\rho_{\iota}$,
\begin{align*}
dm^2 \nu^o \cdot \tilde{\nu}^e &=\sum_{\iota \in I_o} \nu^o_{\iota}
\{\lambda_{\iota} \sum_{\eta \in I_o , \eta \ge \iota} \tilde{\nu}^o_{\eta}
\mu_{\eta}
 +\mu_{\iota}\sum_{\eta \in I_o, \eta <\iota} \tilde{\nu}^o_{\eta}
\lambda_{\eta}\}\\
&=\sum_{\iota \in I_o} \nu^o_{\iota}[-\lambda_{\iota}\sigma (\tilde{\nu}^o)
 + dm\{\rho_{\iota}\sum_{\eta \in I_o, \eta <\iota}\tilde{\nu}^o_{\eta}
\lambda_{\eta}
 +\lambda_{\iota}\sum_{\eta \in I_o, \eta \ge \iota} \tilde{\nu}^o_{\eta}
\rho_{\eta}\}]\\
&=-\sigma (\nu^o)\sigma (\tilde{\nu}^o)
 +dm\sum_{\iota \in I_o} \nu^o_{\iota}
 \{ \rho_{\iota}\sum_{\eta \in I_o, \eta <\iota} \tilde{\nu}^o_{\eta}
\lambda_{\eta}
 +\lambda_{\iota}\sum_{\eta \in I_o, \eta \ge \iota} \tilde{\nu}^o_{\eta}
\rho_{\eta}\}
\end{align*}

By the formula $\rho_{\iota}=-m\bar{\lambda}_{\iota}+a\lambda_{\iota}$, we can
get

\begin{align*}
&\rho_{\iota}\sum_{\eta \in I_o, \eta <\iota}\tilde{\nu}^o_{\eta}\lambda_{\eta}
 +\lambda_{\iota}\sum_{\eta\in I_o, \eta \ge
\iota}\tilde{\nu}^o_{\eta}\rho_{\eta}\\
 =&a\lambda_{\iota}\cdot\sigma(\tilde{\nu}^o)
 -m(\bar{\lambda}_{\iota}\sum_{\eta\in I_o, \eta <\iota}
\tilde{\nu}^o_{\eta}\lambda_{\eta}
 +\lambda_{\iota}\sum_{\eta\in I_o,
\eta\ge\iota}\tilde{\nu}^o_{\eta}\bar{\lambda}_{\eta}
\end{align*}
Hence we get
$$dm^2\nu^o\cdot\tilde{\nu}^o = (dma-1)\sigma(\nu^o)\sigma (\tilde{\nu}^o)
-dm^2 \cal O (\nu^o, \tilde{\nu}^o)$$
Simarlarly we get
\begin{equation*}
dm^2 \nu^e\cdot\tilde{\nu}^e
=-(dma+1)\tau(\nu^e)\tau(\tilde{\nu}^e)
 +dm^2 \cal E (\nu^e, \tilde{\nu}^e)
\end{equation*}
Summarizing all the above formulas, we get the first equality.

 Next for the non-normal case, we can get
$m\nu\cdot\tilde{\nu}=
a\sigma(\nu)\sigma(\tilde{\nu})-m\cal O(\nu,\tilde{\nu})$
for $\nu,\tilde{\nu}\in\cal L_o$,
$m\nu\cdot\tilde{\nu}=
-a\tau(\nu)\tau(\tilde{\nu})+m\cal E(\nu,\tilde{\nu})$
for $\nu,\tilde{\nu}\in\cal L_e$.
We leave the details for the reader.
\end{pf}

Note that $\nu\cdot\tilde{\nu}=\cal E (\nu^e, \tilde{\nu}^e)-\cal O (\nu^o,
\tilde{\nu}^o)$ if $\sigma (\nu^o)=\tau (\nu^e)$ and
$\sigma (\tilde{\nu}^o)=\tau (\tilde{\nu}^o)$ hold.

\begin{cor}\label{corint}
If $\nu^o\le\tilde{\nu}^o, \nu^e\le\tilde{\nu}^e, \sigma (\nu^o)=\tau (\nu^e),
\sigma (\tilde{\nu}^o)=\tau (\tilde{\nu}^e)$ and
$\bar{\sigma}(\tilde{\nu}^o)=\bar{\tau}(\tilde{\nu}^e)$ hold,
then $\nu\cdot\tilde{\nu}=0$.
\end{cor}

\begin{pf}
This can be easily checked by the above lemma.
\end{pf}

\begin{dfn}
For $\iota =(i, j)\in I$ such that $i\not= k+1$,
put $$\varphi (\iota)=-\delta^{i, 1}+\delta^{i, j}+(j-1)\delta^{i+1, j}$$
\end{dfn}

\begin{lem}\label{lemvarphi}
\begin{align*}
&\sigma (\varphi (\iota )^o)=\tau(\varphi (\iota)^e)=(j-1)(P_{i-1}+Q_{i-1})\\
&\bar{\sigma}(\varphi(\iota)^o)=\bar{\tau}(\varphi (\iota)^e)=(j-1)P_{i-1}\\
&\varphi (\iota)^2=j-1
\end{align*}
\end{lem}

\begin{pf}
We can get the first and the second formula by direct calculation.
By the formula $P_{i-2}Q_{i-1}-Q_{i-2}P_{i-1}=(-1)^i$,
we get $$\cal E (\varphi (\iota)^e, \varphi
(\iota)^e)=(j-1)^2P_{i-1}(P_{i-1}+Q_{i-1})$$
and $$\cal O (\varphi (\iota)^o, \varphi (\iota)^o)=(j-1)^2
P_{i-1}(P_{i-1}+Q_{i-1})-j+1$$
Thus by the above corollary, we can get the third formula.
\end{pf}

\begin{dfn}
Let $\iota =(i_1, j_1)$, $\eta=(i_2, j_2)$ be elements in $I$ such that
$i_2\not= k+1$, the parity of $i_1$ coincides the one of $i_2$
and $\iota\le (i_2, 1)$.
For such pair $(\iota, \eta)$, we define $\psi(\iota, \eta)$ as follows
$$\psi (\iota, \eta)=
\begin{cases}
 -\delta^{\iota^l}+\delta{\iota}
 +\sum
 \begin{Sb}
 \iota\le (i, 1)\le\eta \\ i\ \text{\em{odd}}, i\not= 1
 \end{Sb}
 q_{i-1}\delta^{i,1} -\delta^{i_2, 1}+\delta^{\eta}+j_2 \delta^{i_2+1, 1}\quad
&
 (\iota\in I_o)\\
 -\delta^{\iota^r}+\delta{\iota}
 +\sum
 \begin{Sb}
 \iota\le (i, 1)\le\eta \\ i\ \text{\em{even}}
 \end{Sb}
 q_{i-1}\delta^{i,1} -\delta^{i_2, 1}+\delta^{\eta}+j_2 \delta^{i_2+1, 1}\quad
&
  (\iota\in I_e)
\end{cases}
$$
\end{dfn}

\begin{lem}\label{lempsi}
Let $(\iota, \eta =(i'', 1))$ be the pair for which $\psi$ can be defined.
Then
$$\sigma (\psi (\iota, \eta)^o)=\tau (\psi (\iota, \eta
)^e)=P_{i''+1}+Q_{i''+1}$$
$$\bar{\sigma}(\psi (\iota, \eta )^o)=P_{i''+1},\quad
 \bar{\tau}(\psi (\iota, \eta )^e)=
 \begin{cases}
  P_{i''-1}\qquad & (\iota \not= (2, 1))\\
  P_{i''-1}+1\quad & (\iota = (2, 1))
 \end{cases}
$$
$$\psi (\iota, \eta )^2 =
 \begin{cases}
 1+\sum
 \begin{Sb}
 i'+1\le i\le i''-1 \\ i\ \text{\em{even}}
 \end{Sb}
 q_i\qquad & (\iota\in I_o)\\
 1+\sum
 \begin{Sb} i'+1\le i\le i''-1 \\ i\ \text{\em{odd}}
 \end{Sb}
 q_i\qquad & (\iota\in I_e)
 \end{cases}
$$
\end{lem}

\begin{pf}
Since the calculation is the same, we show the outline of it for the case
$\iota \in I_o$.
We can easily calculate $\sigma$, $\tau$, $\bar{\sigma}$, and $\bar{\tau}$.
Hence
$$\psi (\iota, \eta )^2 =\cal E (\psi (\iota, \eta)^e, \psi (\iota, \eta
)^e)-\cal O (\psi (\iota, \eta)^o, ( \psi (\iota, \eta )^o)$$
By the definion,
$$\cal E (\psi (\iota, \eta)^e, \psi (\iota,
\eta)^e)=P_{i''-1}(P_{i''-1}+Q_{i''-1})$$
To calculate $\cal O (\psi (\iota, \eta)^o, \psi (\iota, \eta)^o)$, note that
$$-\lambda_{\iota^l}+\lambda_{\iota}
+\sum
\begin{Sb}
i'+2\le i\le h \\ i\ \text{odd}
\end{Sb}
q_{i-1}\lambda^i_1=P_{h-1}+Q_{h-1}$$
and
$$\sum
\begin{Sb}
h+2\le i\le i''\\ i\ \text{odd}
\end{Sb}
q_{i-1}\bar{\lambda}^i_1=P_{i''-1}-P_{h-1}$$
for odd $h$ such that $i'+2\le h\le i''$. Then
\begin{align*}
\cal O (\psi (\iota, \eta)^o, \psi (\iota,\eta)^o)=
&-\lambda_{\iota ^l}\bar{\sigma}(\psi (\iota,
\eta)^o)+(-\lambda_{\iota^l}+\lambda_{\iota})\bar{\lambda}_{\iota}\\
&+\lambda (\sum
\begin{Sb}
i'+2\le i \le i''\\ i\ \text{odd}
\end{Sb}
q_{i-1}\bar{\lambda}^i_1)\\
&+\sum\begin{Sb} i'+2\le h \le i''\\ h\ \text{odd}
\end{Sb}
q_{h-1}\{ \bar{\lambda}^h_1 (P_{h-1}+Q_{h-1})+\lambda (P_{i''-1}-P_{h-1})\}
\end{align*}

Here
\begin{align*}
&-\lambda_{\iota^l}\bar{\sigma}(\psi (\iota,
\eta)^o)+(-\lambda_{\iota^l}+\lambda_{\iota})\bar{\lambda}_{\iota}
 +\lambda_{\iota} (\sum
\begin{Sb}
i'+2\le i\le i''\\ i\ \text{odd}
\end{Sb}
q_{i-1}\bar{\lambda}^i_1 )\\
=&
(-\lambda_{\iota}+P_{i'-1}+Q_{i'-1})P_{i''-1}+(P_{i'-1}+Q_{i'-1})\bar{\lambda}_{\iota}
 +\lambda_{\iota} (P_{i''-1}-P_{i'-1})\\
=&(P_{i'-1}+Q_{i'-1})P_{i''-1}+(P_{i'-1}+Q_{i'-1})\bar{\lambda}_{\iota}+\lambda_{\iota} P_{i-1}\\
=&(P_{i'-1}+Q_{i'-1})P_{i''-1}+Q_{i'-1}P_{i'-2}-P_{i'-1}Q_{i'-2}\\
=&(P_{i'-1}+Q_{i'-1})P_{i''-1}-1
\end{align*}

and
\begin{align*}
&\bar{\lambda}^h_1 (P_{h-1}+Q_{h-1})+\lambda^h_1 (P_{i'-1}-P_{h-1})\\
=&(P_{h-2}+Q_{h-2})P_{i''-1}-1
\end{align*}

Hence
\begin{align*}
&\cal O (\psi (\iota, \eta)^o, \psi (\iota, \eta)^o)\\
=&(P_{i'-1}+Q_{i'-1})P_{i''-1}-1+\sum
 \begin{Sb} i'+2\le h\le i''\\ h\ \text{odd} \end{Sb}
 q_{h-1}\{(P_{h-2}+Q_{h-2})P_{i''-1}-1\}\\
=&-1-\sum
 \begin{Sb} i'+2\le h \le i''\\ h\ \text{odd}\end{Sb}
 q_{h-1}+P_{i''}\{(P_{i'+1}+Q_{i'-1})+
 \sum\begin{Sb} i'+2\le h\le i''\\ h\ \text{odd}\end{Sb}
q_{h-1}(P_{h-2}+Q_{h-2})\}\\
=&-1-\sum
 \begin{Sb}i'+1\le i\le i''-1\\ i\ \text{even}\end{Sb}
 q_i+P_{i''-1}(P_{i''-1}+Q_{i''-1})
\end{align*}

Thus we have done.
\end{pf}

\begin{cor}\label{corpsi}
Let $(\iota,\eta)$ be a pair for which $\psi$ can be defined. Then
$$\sigma (\psi (\iota, \eta)^o)
 =\tau (\psi (\iota, \eta)^e)
 =j''(P_{i''-1}+Q_{i''-1})$$

\begin{equation*}
 \bar{\sigma}(\psi (\iota, \eta)^o)=
 j''P_{i''-1 },\quad
 \bar{\tau}(\psi (\iota, \eta)^e)=
 \begin{cases}
  j''P_{i''-1}\qquad &(\iota\not= (2, 1))\\
  j''P_{i''-1}+1\quad &(\iota =(2, 1))
 \end{cases}
\end{equation*}

\begin{equation*}
 \psi (\iota, \eta)=
 \begin{cases}
  j''+\sum
 \begin{Sb}i'+1\le i\le i''-1\\ i\ \text{\em{even}}\end{Sb}
 q_i\quad &(if\ \iota \in I_o)\\
  j''+\sum
 \begin{Sb}i'+1\le i\le i''-1\\ i\ \text{\em{odd}}\end{Sb}
 q_i\quad &(if\ \iota \in I_e)
 \end{cases}
\end{equation*}

\end{cor}

\begin{pf}
Note that $\psi (\iota, \eta)=\psi (\iota, (i'', 1))+\varphi (\eta)$.
Hence
$$\sigma (\psi (\iota, \eta)^o)
 =\sigma (\psi (\iota, \eta)^o)+\sigma (\varphi (\eta)^o)
 =j''(P_{i''-1}+Q_{i''-1})$$

$\tau$, $\bar{\sigma}$ and $\bar{\tau}$ are simarlarly calculated.
By Corollary \ref{corint}, $\psi(\iota,(i'', 1))\cdot \varphi(\eta)=0$.
Hence
$$\psi(\iota, \eta)^2=\psi(\iota, (i'', 1))^2+\varphi(\eta)^2$$
Thus the formula follows from Lemma \ref{lemvarphi} and \ref{lempsi}.
\end{pf}

\begin{dfn}
For $(i, j)\in I$ such that $1\le i\le k-1$ and $1\le j\le q_i$,
we put
$$\theta (i, j) =-\delta^{i, q_i-j+1}+\delta^{i+2, 1}
  +j\delta^{i+1, 1}$$
\end{dfn}

\begin{lem}\label{lemtheta}
$$\sigma (\theta (i, j)^o)
  =\tau (\theta (i, j)^e)
  =j(P_{i-1}+Q_{i-1})$$

$$\bar{\sigma}(\theta (i, j)^o)
  =\bar{\tau}(\theta (i, j)^e)
  =jP_{i-1}$$

$$\theta (i, j)=j$$

\end{lem}

\begin{pf}
We will only show the outline of the calculation
for odd $i$ since it is similar for even $i$.
We can easily get the formulas for
$\sigma$, $\tau$, $\bar{\sigma}$ and $\bar{\tau}$.
Thus by Lemma \ref{lemint},
$$\theta (i, j)^2 =\cal E (\theta (i,j)^e, \theta (i, j)^e)
-\cal O (\theta (i, j)^o, \theta (i, j)^o)$$

By the definition,
$$\cal E (\theta (i, j)^e, \theta (i, j)^e)=
 j^2P_{i-1}(P_{i-1}+Q_{i-1})$$
and
\begin{align*}
&\cal O (\theta (i, j)^o, \theta (i, j)^o)\\
=&\lambda^i_{q_i-j+1}(\bar{\lambda}^i_{q_i-j+1}-\bar{\lambda}^{i+2}_1)
  +(-\lambda^i_{q_i-j+1}+\lambda^{i+2}_1)\bar{\lambda}^{i+2}_1\\
=&-jP_{i-1}\lambda^i_{q_i-j+1}
 +j(P_{i-1}+Q_{i-1})\bar{\lambda}^{i+2}_1\\
=&j^2P_{i-1}(P_{i-1}+Q_{i-1})+j(P_iQ_{i-1}-Q_iP_{i-1})\\
=&j^2P_{i-1}(P_{i-1}+Q_{i-1})-j
\end{align*}

Hence we have done.
\end{pf}

For a positive integer $n$ which is smaller than $m-(P_{k-1}+Q_{k-1})$,
we put
$$i(n)=\max\{i|0\le i \le k-1, P_i+Q_i\le n\}$$

\begin{prop}\label{propnu1}
If $0<n<m-(P_{k-1}+Q_{k-1})$,
$$\nu (n)^2\le \dfrac{n}{P_{i(n)}+Q_{i(n)}}$$
\end{prop}

\begin{pf}
We use the induction on $i(n)$.
If $i(n)=0$, it can be easily checked that $\nu (n)^2 =n$.
Let $i$ be an integer such that $1\le i\le k-1$.
Assume that the ineqality holds for all $n$ such that $i(n)<i$.
We will show that the inequality holds for $n$
such that $i(n)=i$ under this assumption.
We also assume $i$ is odd since the proof is similar for even $i$.

Write $n=j(P_i+Q_i)+n'$ such that $0\le n'<P_i+Q_i$.
If $n'=0$, we can check
 (by the definition of $\lambda$- and $\mu$- expansion)
 that $$\nu (n)=\psi ((2, 1), (i+1, j))$$
Thus by Corollary \ref{corpsi},
$$\nu (n)^2\le j$$
Hence we have done in this case.
Thus we assume $n'>0$.
We divide the proof into two cases as follows
(i)$n'<P_{i-1}+Q_{i-1}$ (ii)$n'\ge P_{i-1}+Q_{i-1}$

(i) Put
$$\eta =\min\{\iota\in\bar{I}_e|\iota '<\iota\ for\ all\ \iota '\in \bar{I}_e
 \ such\ that\ \nu (n')_{\iota '}\not= 0\}$$
We can check (by the definition of $\lambda$- and $\mu$- expantion)
$$\nu (n)=\nu (n')+\psi (\eta, (i+1, j))$$
Since $\nu (n')^o\le \psi (\eta, (i+1, j))^o$
and $\nu (n')^e\le \psi (\eta, (i+1, j))^e$ hold,
thus by Corollary \ref{corint},
$$\nu (n)^2=\nu (n')^2+\psi (\eta, (i+1, j))^2$$
By Corollary \ref{corpsi} and the induction hypothesis,
\begin{align*}
 \nu (n)^2(P_{i}+Q_{i})&\ge \nu (n')^2(P_{i}+Q_{i})
                               +j(P_{i}+Q_{i})\\
 &\ge \nu (n')^2(P_{i(n')}+Q_{i(n')})+j(P_{i}+Q_{i})\\
 &\ge n'+j(P_i+Q_i)\\
 &=n
\end{align*}
Thus we have done.
\par
(ii) We can check
$$\nu (n)=\nu (n')+\varphi (i+1, j+1)$$
$$\nu (n')^o\le \varphi (i+1, j+1)^o,\quad
  \nu (n')^e\le \varphi (i+1, j+1)^e$$

Thus we can get the inequality by the similar way to (i)
using Lemma \ref{lemvarphi} and the induction hypothesis,
\end{pf}

For $n$ such that $m-(P_{k-1}+Q_{k-1})\le n\le m-1$,
we define $i(n)$ and $j(n)$ as follows
$$i(n)=\min\{i|0\le i\le k-1, m-n\le P_i+Q_i\}$$
$$j(n)=\lceil \dfrac{m-n}{P_{i(n)-1}+Q_{i(n)-1}} \rceil -1$$

\begin{lem}\label{lemnutheta}
Let $n$, $i$, $j$ be positive integers such that
$m-(P_{k-1}+Q_{k-1})\le n\le m-1$,
$i\le k-1$, $j\le q_i$ and $i(n)\le i-1$.
Then $\nu (n)\cdot\theta (i, j)=j$.
\end{lem}

\begin{pf}
We only show the proof for even $i$.
By Lemma \ref{lemint},
$$\nu (n)\cdot\theta (i, j)
 =\cal E (\nu (n)^e, \theta (i, j)^e)
  -\cal O (\nu (n)^o, \theta (i, j)^o)$$

First we calculate $\cal E (\nu (n)^e, \theta (i, j)^e)$.
Note that $m-(P_{i-1}+Q_{i-1})\le n\le m-1$.
Thus the $\mu$-expantion of $n$ is as follows

\begin{equation*}
n=\sum\begin{Sb} 2\le h\le i-2\\ h\ \text{even}\end{Sb}
  (l_{h-2}\mu ^h_1 +\mu ^h_{j_h})
   +l_{i-1}\mu ^i_1
  +\sum\begin{Sb}i+2\le h\le k-2\\ h\ \text{even}\end{Sb}
   q_{h-1}\mu ^h_1
   +q_{k-1}\mu ^k_1 +\mu ^k_{q_k}
\end{equation*}

Since $\sum\begin{Sb}\iota\in \bar{I}_e\\ \iota\ge (i+2, 1)\end{Sb}
    \nu (n)^e_{\iota}\mu_{\iota}
   =m-(P_{i-1}+Q_{i-1})$,
we get
$$\sum\begin{Sb} \iota\in \bar{I}_e\\ \iota\le (i, 1)\end{Sb}
  \nu (n)^e_{\iota} \mu_{\iota}
   =n-m+P_{i-1}+Q_{i-1}$$

Using this formula and the formula
$\sum_{\iota\in \bar{I}_e, \iota\ge (i+2, 1)}
 \nu (n)^e_{\iota}\bar{\mu}_{\iota}
 =a-P_{i-1}$,
we get

\begin{align*}
&\cal E (\nu (n)^e, \theta (i, j)^e)\\
=&\cal E (\sum_{\iota\in\bar{I}_e, \iota\le (i,1)}
\nu (n)^e_{\iota}\delta^{\iota}
+\sum_{\iota\in I_e, \iota\ge (i+2, 1)}
\nu (n)^e_{\iota}\delta ^{\iota},
\theta (i, j)^e)\\
=&\tau (\sum _{\iota\in\bar{I}_e, \iota\le (i, 1)}
\nu (n)^e_{\iota}\delta ^{\iota}
\cdot\bar{\tau}(\theta (i, j)^e)
+\bar{\tau}(\sum_{\iota\in\bar{I}_e, \iota\le (i+2, 1)}
\nu (n)^e_{\iota}\delta^{\iota}
\cdot\tau (\theta (i, j)^e)\\
=&(n-m+P_{i-1}+Q_{i-1})\cdot jP_{i-1}
+(a-P_{i-1})\cdot j(P_{i-1}+Q_{i-1})\\
=&j(n-m)P_{i-1}+ja(P_{i-1}+Q_{i-1})
\end{align*}

Next we calculate $\cal O (\nu (n)^o, \theta (i, j)^o)$.
The $\lambda$-expansion of $n$ is as follows

\begin{equation*}
n=\lambda^1_{j_1}+\sum\begin{Sb} 3\le h\le i-1 \\ h\ \text{odd} \end{Sb}
    l_{h-1}\lambda^h_1+\lambda^h_{j_h}+l_{i+1}\lambda^{i+1}_1
   +\sum\begin{Sb} i+3\le h\le k-1 \\ h\ \text{odd} \end{Sb}
    q_{h-1}\lambda^h_1+q_k\lambda^{k+1}_1
\end{equation*}

Since $\sum\begin{Sb} \iota\in I_o \\ \iota>(i+1, 1)\end{Sb}
 \nu (n)^o_{\iota}\lambda_{\iota}
 =m-(P_{i-1}+Q_{i-1})$,
we get
$$\sum\begin{Sb} \iota\in I_o \\ \iota\le (i+1,1)\end{Sb}
 \nu (n)^o_{\iota}\lambda_{\iota}=
 n-m+P_i+Q_i$$
Using this formula and
$\sum\begin{Sb}\iota\in I_o\\ \iota >(i+1, 1)\end{Sb}
\nu (n)^o_{\iota}\bar{\lambda}_{\iota}
=a-P_{i-1}$,
we can get

\begin{align*}
&\cal O (\nu (n)^o, \theta (i, j)^o)\\
=&\cal O (\sum
\begin{Sb} \iota\in I_o\\ \iota\le (i+1, 1)\end{Sb}
 \nu (n)^o_{\iota}\delta^{\iota}, \theta (i, j)^o )
 +\cal O (\sum
 \begin{Sb}\iota\in I_o\\ \iota>(i+1, 1)\end{Sb}
 \nu (n)^o_{\iota}\delta^{\iota}, \theta (i, j)^o)\\
=&\sigma (\sum
 \begin{Sb} \iota\in I_o\\ \iota\le (i+1, 1)\end{Sb}
 \nu (n)^o_{\iota}\delta^{\iota})
 \bar{\sigma}(\theta (i, j)^o)
 +\bar{\sigma}(\sum
 \begin{Sb} \iota\in I_o\\ \iota>(i+1, 1)\end{Sb}
 \nu (n)^o_{\iota}\delta^{\iota})
 \sigma (\theta (i, j)^o)\\
=&(n-m+P_i+Q_i)\cdot jP_{i-1}+(a-P_i)\cdot j(P_{i-1}+Q_{i-1})\\
=&j(n-m)P_{i-1}+ja(P_{i-1}+Q_{i-1})
  -j(P_iQ_{i-1}-Q_iP_{i-1})\\
=&j(n-m)P_{i-1}+ja(P_{i-1}+Q_{i-1})-j
\end{align*}

Thus we have done.
\end{pf}

\begin{lem}
$$\nu (m-1)^2=\sum_{1\le h\le k}q_h$$
\end{lem}

\begin{pf}
Note that $\lambda$- and $\mu$- expansion of $m-1$ is as follows

\begin{align*}
m-1&=\lambda^1_0+\sum
 \begin{Sb}3\le h\le k-1\\ h\ \text{odd}\end{Sb}
(q_{h-1}\lambda^h_1+\lambda^h_0)+q_k\lambda^{k+1}_1\\
&=\sum
 \begin{Sb}2\le h \le k-2\\ h\ \text{even}\end{Sb}
 (q_{h-1}\mu^h_1+\mu^h_0)+q_{k-1}\mu^k_1+\mu^k_{q_k}
\end{align*}

We leave the rest of calculation for the reader's exercise.
\end{pf}

\begin{prop}\label{propnu2}
Let $n$ be an integer such that
$m-(P_{k-1}+Q_{k-1})\le n\le m-1$.
Then
$$\nu (n)^2\ge \sum_{i(n)\le h\le k}q_h-j(n)$$
\end{prop}

\begin{pf}
We use the induction on $i(n)$.
If $i(n)=0$, then $n=m-1$, Thus we have already done in the above lemma.
Let $i$ be a positive integer and assume that the inequality holds
for $n'$ such that $i(n')<i$.
Let $n$ be an integer such that the inequality holds for this $n$.
Put $n'=n+j(P_{i-1}+Q_{i-1})$.
Then $i(n')<i$.
We can check
$$\nu (n)=\nu (n')-\theta (i, j)$$

Thus by Lemma \ref{lemtheta} and \ref{lemnutheta}, we can get
\begin{align*}
\nu (n)^2&=\nu (n')^2-2\nu (n')\cdot \theta (i, j)
  +\theta (i, j)^2\\
 &=\nu (n')^2-j
\end{align*}

By the induction hypothesis,
$$\nu (n')\ge \sum_{i(n')+1\le h\le k}q_h$$
Thus we have done.
\end{pf}

(Proof of the Theorem \ref{thmbound})\\
{}From Proposition \ref{propnu1} and Proposition \ref{propnu2},
we know the thorem holds if $D$ and $D'$ is general members in $|nK_X|$.
Thus by Corollary \ref{corgenint}, we have proved the thorem.




\begin{thebibliography}{[K-SB]}
\bibitem[ESN]{esn}
H. Esnault,
{\it Reflexive modules on quotient surface singularities},
J. Reine Angew. Math. 362 (1985), 63--71


\bibitem[K-SB]{ksb}
J. Koll\'{a}r and N. Shepherd-Barron,
{\it Threefolds and deformations of surface singularities},
Invent. Math. 91 (1988), 299--338
\end{thebibliography}
\end{document}